\documentclass[reqno, a4paper]{amsart}
\usepackage{amsmath}
\usepackage{amssymb}
\usepackage{amsthm}

\usepackage[scale=0.9]{geometry}

\usepackage{natbib}
\usepackage{bibentry} 

\usepackage[english]{babel}
\usepackage[utf8]{inputenc}

\usepackage{subfig}
\usepackage{graphicx}

\usepackage[unicode]{hyperref}
\usepackage[active]{srcltx} 

\usepackage{mathbbol}
\usepackage{bm} 
\usepackage{stmaryrd} 
\usepackage{MnSymbol} 

\usepackage{units}
\usepackage{tensor}
\usepackage{accents}

\usepackage{enumitem}

\usepackage{lineno}

\usepackage{comment}

\usepackage{placeins}

\usepackage{todonotes}

\usepackage{algorithm}
\usepackage{algpseudocode}

\usepackage{fvextra} 

\newcommand{\bigo}[1]{\ensuremath{O\left(#1 \right)}}

\newcommand{\bydefinition}{\mathrm{def}}

\newcommand{\diff}{\mathrm{d}}

\renewcommand{\vec}[1]{\ensuremath{\mathbf{#1}}}

\makeatletter
\@ifpackageloaded{bm}%
{\renewcommand{\vec}[1]{\ensuremath{\bm{#1}}}%
}{%
\relax
}
\makeatother
\newcommand{\tensorq}[1]{\ensuremath{\mathbb{#1}}}      
\newcommand{\tensorc}[1]{\ensuremath{\mathrm{#1}}}      

\newcommand{\transpose}[1]{#1^\top}

\newcommand{\identity}{\ensuremath{\tensorq{I}}} 

\makeatletter
\@ifpackageloaded{bm}%
{%
}{%

}

\@ifpackageloaded{bm}%
{%
 
}{%

}

\@ifpackageloaded{bm}%
{%
}{%

}
\makeatother

\newcommand{\kdelta}[1]{\tensor{\delta}{#1}}

\newcommand{\diracdelta}{\delta}

\newcommand{\R}{\ensuremath{{\mathbb R}}}

\newcommand{\N}{\ensuremath{{\mathbb N}}}

\newcommand{\dd}[2]{\ensuremath{\frac{\diff {#1}}{\diff {#2}}}}

\newcommand{\ddd}[2]{\ensuremath{\frac{\diff^2 {#1}}{\diff {#2}^2}}}

\makeatletter
\@ifpackageloaded{tensor}
{

}{%

}
\makeatother

\makeatletter
\@ifpackageloaded{tensor}
{

}{%

}
\makeatother

\makeatletter
\@ifundefined{volume}{%
}%
{%
}
\makeatother

\makeatletter

\@ifpackageloaded{MnSymbol} 
{
 
}{%
 
}

\@ifpackageloaded{MnSymbol} 
{
   
}{%
   
}
\makeatother

\newcommand{\vectordot}[2]{\ensuremath{#1 \bullet #2}}

\newcommand{\tensorschur}[2]{\ensuremath{#1 \circ #2}} 

\newcommand{\vectordotalt}[3]{\ensuremath{#1 \bullet_{#3} #2}} 

\newcommand{\fascalalt}[2]{\left\langle #1, #2 \right\rangle} 

\newcommand{\DMatrixc}[1]{\tensor{\tensorc{D}}{_{#1}}} 
\newcommand{\DMatrix}[1]{\tensorq{D}_{#1}} 

\newcommand{\DDMatrixBC}[1]{\tensorq{D}_{#1, \, \mathrm{bc}}^2} 

\newcommand{\DDMatrixTilde}[1]{\widetilde{\tensorq{D}_{#1}^2}} 

\newcommand{\GMatrixc}[1]{\tensor{\tensorc{G}}{_{#1}}} 
\newcommand{\GMatrix}[1]{\tensorq{G}_{#1}} 

\newcommand{\GMatrixBC}[1]{\tensorq{G}_{#1, \, \mathrm{bc}}} 

\newcommand{\GFunction}{g} 
\newcommand{\GOperator}{\mathcal G} 

\newcommand{\bweight}{\lambda} 
\newcommand{\qweight}{w} 

\newcommand{\ChebyshevT}[1]{\mathrm{T}_{#1}} 

\newcommand{\vecelem}[1]{#1} 

\DeclareMathOperator{\DCTI}{DCT-I}
\DeclareMathOperator{\INTEGRATE}{INTEGRATE}
\DeclareMathOperator{\EXTEND}{EXTEND}
\DeclareMathOperator{\REDUCE}{REDUCE}

\title[Fast construction of discrete Green operator]{Fast construction of the discrete Green operator for a second order ordinary differential equation}

\author{Jan Blechta}
\address{
Faculty of Mathematics and Physics\\
Charles University\\
Sokolovsk\'a 83\\
Praha 8 -- Karl\'{\i}n\\
CZ 186\;75\\
Czech Republic
}
\email{blechta@kalrin.mff.cuni.cz}

\author{V\'{\i}t Pr\r{u}\v{s}a}
\date{\today}
\address{
Faculty of Mathematics and Physics\\
Charles University\\
Sokolovsk\'a 83\\
Praha 8 -- Karl\'{\i}n\\
CZ 186\;75\\
Czech Republic
}
\email{prusv@karlin.mff.cuni.cz}

\author{Ladislav Trnka}
\address{
Faculty of Mathematics and Physics\\
Charles University\\
Sokolovsk\'a 83\\
Praha 8 -- Karl\'{\i}n\\
CZ 186\;75\\
Czech Republic
}
\email{ladtrnjr@volny.cz}

\author{Karel T\r{u}ma}
\address{
Faculty of Mathematics and Physics\\
Charles University\\
Sokolovsk\'a 83\\
Praha 8 -- Karl\'{\i}n\\
CZ 186\;75\\
Czech Republic
}
\email{ktuma@karlin.mff.cuni.cz}

\thanks{V\'{\i}t Pr\r{u}\v{s}a thanks the Czech Science Foundation, grant number 20-11027X, for its support.}

\numberwithin{equation}{section}

\begin{document}

\begin{abstract}
  We consider linear second order differential equation $\ddd{y}{x} = f$ with zero Dirichlet boundary conditions. At the continuous level this problem is solvable using the Green function, and this technique has a counterpart on the discrete level. The discrete solution is represented via an application of a matrix---the Green matrix---to the discretised right-hand side, and we propose an algorithm for fast construction of the Green matrix. In particular, we discretise the original problem using the spectral collocation method based on the Chebyshev--Gauss--Lobatto points, and using the discrete cosine transformation we show that the corresponding Green matrix is fast to construct even for large number of collocation points/high polynomial degree. Furthermore, we show that the action of the discrete solution operator (Green matrix) to the corresponding right-hand side can be implemented in a matrix-free fashion.



\end{abstract}

\maketitle

\tableofcontents


\section{Introduction}
\label{sec:introduction}
We consider a simple boundary value problem for a second order ordinary differential equation with zero Dirichlet boundary conditions. The task is to find a function $y$ that solves the problem
\begin{subequations}
  \label{eq:1}
  \begin{align}
    \label{eq:2}
    \ddd{y}{x} &= f,  \\
    \label{eq:3}
    \left. y \right|_{x=-1} &= 0, \\
    \label{eq:4}
    \left. y \right|_{x=1} &= 0,
  \end{align}
\end{subequations}
where the function $f$ on the right-hand side of~\eqref{eq:2} is a given function. The predominant approach to the numerical solution of this problem is based on the discretisation of the differential operator on the left-hand side, and on subsequent solution of the system of linear algebraic equations
\begin{equation}
  \label{eq:5}
  \DDMatrixTilde{N} \widetilde{\vec{y}}= \widetilde{\vec{f}},
\end{equation}
where $\DDMatrixTilde{N}$, $\widetilde{\vec{y}}$ and $\widetilde{\vec{f}}$ are discrete counterparts of $\ddd{}{x}$, $y$ and $f$. The discrete problem~\eqref{eq:16} is typically solved using algorithms for solution of systems of linear algebraic equations, which means that the problem is \emph{not} interpreted as the problem of finding the \emph{inverse} to the matrix $\DDMatrixTilde{N}$; one simply wants, for given $\widetilde{\vec{f}}$, to find the corresponding $\widetilde{\vec{y}}$ without explicitly constructing the inverse matrix/the discrete solution operator.

On the continuous level one can, however, explicitly construct the solution operator to~\eqref{eq:1}. The solution operator is constructed using the Green function, see, for example, \cite{lanczos.c:linear}. The Green function $\GFunction\left(x, \xi \right)$ is a function that for a fixed $\xi \in [-1,1]$ solves the problem
\begin{subequations}
  \label{eq:6}
  \begin{align}
    \label{eq:7}
    \ddd{\GFunction\left(x, \xi \right)}{x} &= \diracdelta(x - \xi),  \\
    \label{eq:8}
    \left. \GFunction\left(x, \xi \right) \right|_{x=-1} &= 0, \\
    \label{eq:9}
    \left. \GFunction\left(x, \xi \right) \right|_{x=1} &= 0,
  \end{align}
\end{subequations}
where $\diracdelta$ denotes the Dirac delta distribution. The solution to the problem~\eqref{eq:1} is then obtained by the integration
\begin{equation}
  \label{eq:10}
  y(x) = \int_{\xi=-1}^1 {\GFunction(x,\xi)} f(\xi) \, \diff \xi.
\end{equation}
For further reference we denote this operator as $\GOperator$, that is for function $h$ we define
\begin{equation}
  \label{eq:111}
  \GOperator(h) = _{\bydefinition}  \int_{\xi=-1}^1 {\GFunction(x,\xi)} h(\xi) \, \diff \xi,
\end{equation}
and we refer to this operator as the Green operator/solution operator. Formula~\eqref{eq:10} can be again discretised, and then it admits the representation as the matrix vector multiplication
\begin{equation}
  \label{eq:11}
  \vec{y} = \GMatrix{N} \vec{f},
\end{equation}
where $\GMatrix{N}$, $\vec{y}$ and $\vec{f}$ are discrete counterparts of the integral operator~\eqref{eq:10} and functions $y$ and $f$. The matrix $\GMatrix{N}$, which we refer to as the \emph{Green matrix}, then represents---in the sense discussed in Section~\ref{sec:green-matrix-as}---the ``inverse'' of the matrix $\DDMatrixTilde{N}$.

In what follows we discretise the original problem~\eqref{eq:1} using the spectral collocation method, see, for example, \cite{trefethen.ln:spectral} or~\cite{weideman.ja.reddy.sc:matlab}, and \emph{we propose an algorithm for fast and stable numerical construction of the corresponding Green matrix} $\GMatrix{N}$. This can be seen as an effort complementary to the effort invested in the search for the fast and stable construction of the second order spectral differentiation matrix~$\DMatrix{N}^2$, see~\cite{weideman.ja.reddy.sc:matlab}. The proposed algorithm is implemented in Matlab and Wolfram Mathematica. The former implementation provides interoperability with the \texttt{chebfun} package, which is a standard Matlab package used in, among others, implementation of the spectral collocation method. The later implementation provides interoperability with symbolic analytical tools available in Wolfram Mathematica. Furthermore, in Matlab implementation we also show that the \emph{action of the discrete Green operator/solution operator can be implemented in a matrix-free fashion}.

\section{Fundamentals of spectral collocation method}
\label{sec:fund-spectr-coll}

\subsection{Standard approach}
\label{sec:standard-approach}
The simple boundary value problem~\eqref{eq:1}~is straightforward to solve numerically using the spectral collocation method, see, for example, \cite{boyd.jp:chebyshev}, \cite{trefethen.ln:spectral}, \cite{weideman.ja.reddy.sc:matlab} and~\cite{aurentz.jl.trefethen.ln:block} for comprehensive treatises on the subject matter and for computer codes.

In the spectral collocation method the derivative on the left-hand side is on the discrete level replaced by the second order differentiation matrix, and the equation is enforced at the chosen set of collocation points/nodes. Typically, Chebyshev--Gauss--Lobatto points $\left\{ x_j \right\}_{j=0}^N$,
\begin{equation}
  \label{eq:12}
  x_j =_{\bydefinition} \cos \left( \frac{j \pi}{N}  \right),
\end{equation}
are used for this purpose\footnote{Through the text we are using conventions that are common in the literature on the polynomial interpolation theory. In the actual implementation one typically uses the index set $j=1, \dots, N+1$, which is more suitable for implementation in programming languages that index arrays from one. Matlab and Wolfram Mathematica that we use for proof-of-concept implementation are examples thereof. Furthermore, software packages such as \texttt{chebfun} in fact use the points $x_j =_{\bydefinition} - \cos \left( \frac{(j-1) \pi}{N-1}  \right)$, $j=1, \dots, N+1$, which preserves the natural ordering of the nodes in the sense that $x_0=-1$ and $x_{N+1}=1$.}. In this case the solution is sought in the space of polynomials of degree at most~$N$, and the second order differentiation matrix is given by an explicit formula. The formula for the first order differentiation matrix~$\DMatrix{N}$---the discrete counterpart of $\dd{}{x}$ operator---is known. The elements of the $(N+1) \times (N+1)$ spectral differentiation matrix $\DMatrix{N}$ are given by the formulae
\begin{subequations}
  \label{eq:13}
  \begin{align}
    \label{eq:14}
    \DMatrixc{kj}
    &=
      \frac{\bweight_j}{\bweight_k}
      \frac{1}{x_k-x_j}
      ,
      \qquad j, k = 0, \dots, N, \ j \not = k,
    \\
    \DMatrixc{kk}
    &=
      -
      \sum_{\substack{j=0 \\ j \not = k}}^N
    \DMatrixc{kj}
    ,
    \qquad k=0, \dots, N,
  \end{align}
\end{subequations}
where
\begin{equation}
  \label{eq:15}
  \bweight_j =_{\bydefinition}
  \frac{1}{\prod_{\substack {k=0 \\ k \not = j}}^N \left(x_j-x_k\right)},
\end{equation}
$j=0, \dots, N$ denote the \emph{barycentric weights}, see~\cite{trefethen.ln:spectral,trefethen.ln:approximation}. Subsequently, the matrix $\DMatrix{N,2}$ representing the second order derivative $\ddd{}{x}$ is, in principle, obtained by squaring the spectral differentiation matrix for the first derivative, $\DMatrix{N,2} = \DMatrix{N}^2$. Using various optimisations, see \cite{weideman.ja.reddy.sc:matlab} and the discussion therein, it is possible to compute the spectral differentiation matrices for the first and the second derivative even for large values of $N$.  We recall that in the case of spectral collocation methods the second order differentiation matrix is \emph{dense} and it is \emph{not symmetric}. Instead of \emph{symmetry} the matrix has the \emph{centrosymmetry} property; this lack of symmetry might be puzzling, but it is actually fine, see Appendix~\ref{sec:second-order-spectr}.

The arising system of linear algebraic equations that represents the discrete variant of~\eqref{eq:1} then reads
\begin{equation}
  \label{eq:16}
  \DDMatrixTilde{N} \widetilde{\vec{y}}= \widetilde{\vec{f}},
\end{equation}
where the $(N-1) \times (N-1)$ matrix $\DDMatrixTilde{N}$ denotes the second order spectral differentiation matrix~$\DMatrix{N}^2$ with its first and last rows and columns removed, and the vectors $\widetilde{\vec{y}}$ and $\widetilde{\vec{f}}$ contain values of functions $y$ and $f$ at the inner collocation points $\left\{x_i \right\}_{i=1}^{N-1}$, see~\cite{trefethen.ln:spectral} for details. The system of algebraic equations~\eqref{eq:16} is then typically solved using a direct solver, see, for example~\cite{trefethen.ln:spectral} and the popular \texttt{chebfun} package, which uses Matlab and \texttt{mldivide} function (backslash operator).

\subsection{Green function based approach}
\label{sec:green-function-based}
\begin{figure}[h]
  \centering
  \includegraphics[width=0.45\textwidth]{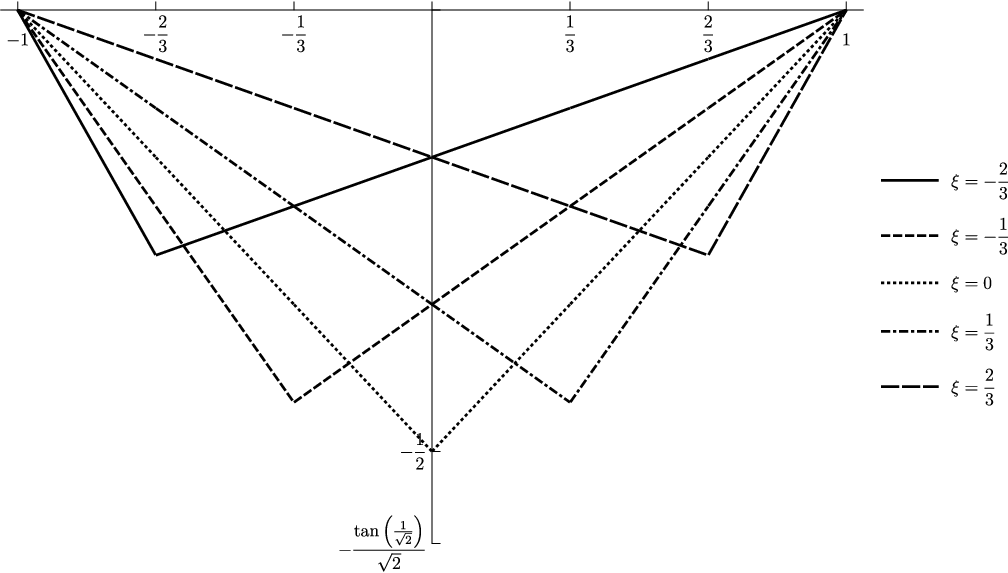}
  \caption{Green function $\GFunction(x, \xi)$ for problem $\ddd{y}{x} = f$; interval $\left(-1, 1\right)$, zero Dirichlet boundary conditions.}
  \label{fig:green-functions}
\end{figure}

The fact that the formula~\eqref{eq:10} gives the solution to the boundary value problem~\eqref{eq:1} is based on the formal manipulation
\begin{equation}
  \label{eq:17}
  \ddd{y}{x}(x)
  =
  \int_{\xi=-1}^1 {\ddd{\GFunction(x,\xi)}{x}} f(\xi) \, \diff \xi
  =
  \int_{\xi=-1}^1 { \diracdelta(x - \xi)} f(\xi) \, \diff \xi
  =
  f(x)
  .
\end{equation}
If necessary, this manipulation can be made rigorous using the standard theory of distributions, see~\cite{schwartz.l:theorie}. The Green function for the problem~\eqref{eq:1} is straightforward to find, and it is given by the formula
\begin{equation}
  \label{eq:18}
  \GFunction(x, \xi)
  =
  \begin{cases}
    \frac{1}{2}(x+1)(\xi -1), &  -1 \leq x \leq \xi  \\
    \frac{1}{2}(x-1)(\xi +1), & \xi < x \leq 1.
  \end{cases}
\end{equation}
Note that the Green function is a \emph{continuous function with a jump discontinuity in the first derivative}. The function is shown in Figure~\ref{fig:green-functions}.

If we want to evaluate the approximate solution $y$ at a collocation point $x_k$, we can use the integral~\eqref{eq:10} and the approximation of the right-hand side $f$ by its polynomial interpolant at the collocation points,
\begin{equation}
  \label{eq:19}
  f(x)
  \approx
  \sum_{i=0}^N f(x_i) l_i(x),
\end{equation}
where $\left\{ l_i \right\}_{i=0}^N$ denotes the Lagrange basis,
\begin{equation}
  \label{eq:20}
  l_j(x)
  =_{\bydefinition}
  \frac{\prod_{\substack{k=0 \\ k \not = j}}^{N} \left(x-x_k\right)}{\prod_{\substack {k=0 \\ k \not = j}}^N \left(x_j-x_k\right)}
\end{equation}
for  $j=0, \dots, N$. Using the polynomial approximation for $f$ in~\eqref{eq:10} we get
\begin{equation}
  \label{eq:22}
  y(x_k)
  =
  \int_{\xi=-1}^1 {\GFunction(x_k,\xi)} f(\xi) \, \diff \xi
  \approx
  \int_{\xi=-1}^1 {\GFunction(x_k,\xi)} \sum_{i=1}^N f(x_i) l_i(\xi) \, \diff \xi
  =
  \sum_{i=0}^N
  \underbrace{
    \left(\int_{\xi=-1}^1 {\GFunction(x_k,\xi)}
      l_i(\xi) \, \diff \xi \right)
  }_{\GMatrixc{ki}}
  f(x_i)
  =
  \sum_{i=0}^N
  \GMatrixc{ki}
  f(x_i)
  .
\end{equation}
Clearly, the application of the Green operator on the right-hand side $f$ corresponds on the discrete level to matrix-vector multiplication. We note that if the right-hand side $f$ is a polynomial of degree at most~$N$, that is if $f$ is identical to its interpolating polynomial, and consequently~\eqref{eq:22} is \emph{exact as well}.

In other words, if we manage to find exact and easy to evaluate formulae for the elements of the discrete \emph{Green matrix}~$\GMatrix{N}$,
\begin{equation}
  \label{eq:23}
  \GMatrixc{ki}
  =_{\bydefinition}
  \int_{\xi=-1}^1
  \GFunction(x_k,\xi)
  l_i(\xi)
  \, \diff \xi
\end{equation}
we obtain a discrete representation of the solution operator to the problem~\eqref{eq:1},
\begin{equation}
  \label{eq:24}
  \vec{y} = \GMatrix{N} \vec{f},
\end{equation}
where $\vec{y}$ and $\vec{f}$ are vectors of pointwise values of functions $y$ and $f$ at \emph{all} collocation points~$\left\{x_i\right\}_{i=0}^N$. Furthermore the solution will be exact provided that the right-hand side is a polynomial of degree at most~$N$. We emphasise that unlike in~\eqref{eq:16} \emph{we do not need to solve a linear system in order to get the approximate discrete solution for the given right-hand side} $f$. All that needs to be done is to \emph{apply} the Green matrix/discrete solution operator to the discrete counterpart of $f$. In fact, it turns out that the action of the discrete solution operator can be performed in a matrix-free fashion.

We now show that the Green matrix can be computed in an easy and numerically feasible way even for large~$N$, that is for large number of collocation points/high polynomial degree. We reiterate that we need an \emph{exact} formula for the integral in~\eqref{eq:23}. Indeed, it would be tempting to use a numerical quadrature, for example the Clenshaw--Curtis quadrature at the given set of collocation points, but this is not an option. The Green function is a function with a jump discontinuity in the first derivative, hence the product $ \GFunction(x_k,\xi) l_i(\xi)$ has bad properties from the perspective of the numerical quadrature, see~\cite{trefethen.ln:exactness} for a careful discussion of numerical quadrature.

Similarly, it would be tempting to use symbolic integration in~\eqref{eq:23}, which is again in principle possible to do explicitly. (The integrand is a piecewise polynomial function, hence the primitive function can be in principle found computationally on-the-fly using software tools for symbolic manipulations.) But such an approach would be prohibitively slow for large $N$, that is for large number of collocation points/high polynomial degree. 

\section{Green matrix construction}
\label{sec:green-matr-constr}
\begin{figure}[h]
  \subfloat[Function $l_i$. \label{fig:discrete-dirac-a}]{\includegraphics[width=0.48\textwidth]{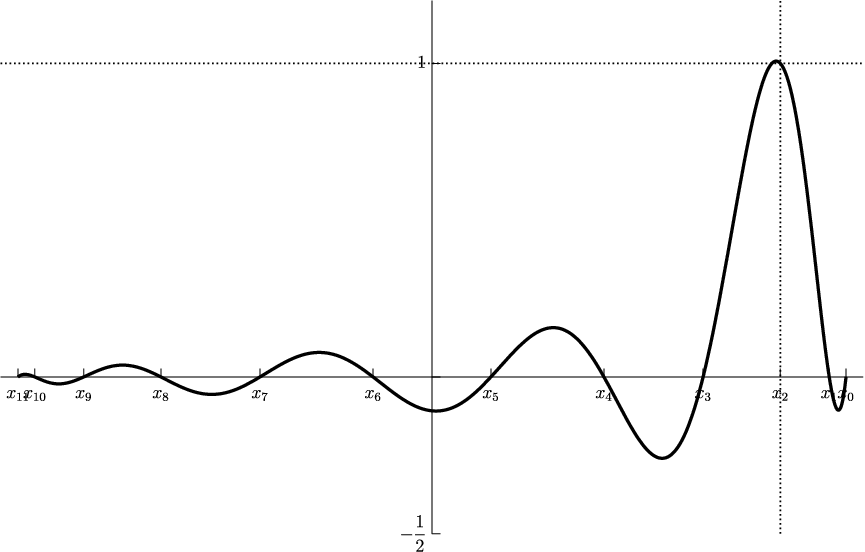}}
  \qquad
  \subfloat[Function $xl_i$. \label{fig:discrete-dirac-b}]{\includegraphics[width=0.48\textwidth]{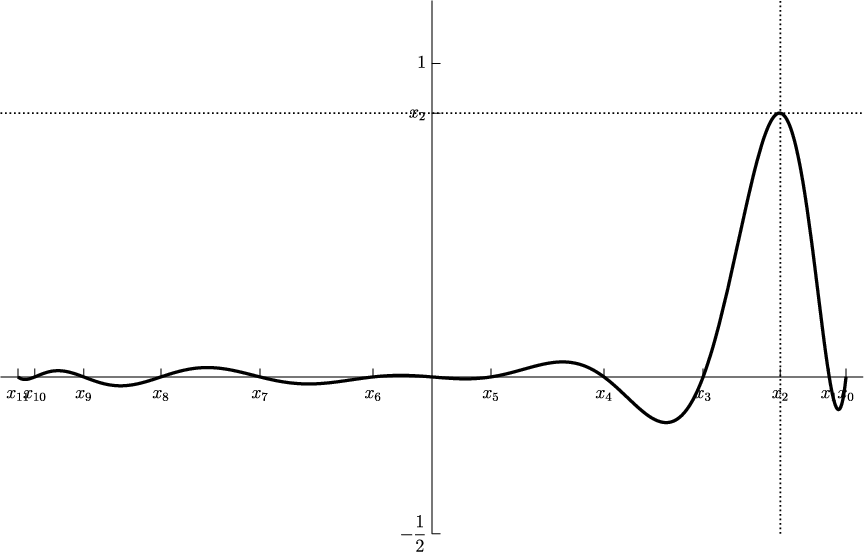}}
  \caption{Integrands in~\eqref{eq:26}, $i=2$, $N=11$.}
  \label{fig:discrete-dirac}
\end{figure}

Using the known analytical formula for the Green function, see~\eqref{eq:18}, we get the elements of the Green matrix~\eqref{eq:23} via the integration 

\begin{multline}
  \label{eq:25}
  \GMatrixc{ki}
  =
  \int_{\xi=-1}^1
  \GFunction(x_k,\xi)
  l_i(\xi)
  \, \diff \xi
  =
  \int_{\xi=x_k}^1 \frac{1}{2}(x_k+1)(\xi -1) l_i(\xi) \, \diff \xi
  +
  \int_{\xi=-1}^{x_k} \frac{1}{2}(x_k-1)(\xi +1) l_i(\xi) \, \diff \xi
  \\
  =
  \frac{1}{2}
  \left(
    x_k+1
  \right)
  \left[
    \int_{\xi = x_k}^1  \xi l_i(\xi) \, \diff \xi
    -
    \int_{\xi = x_k}^1  l_i(\xi) \, \diff \xi
  \right]
  +
  \frac{1}{2}
  \left(
    x_k-1
  \right)
  \left[
    \int_{\xi = -1}^{x_k}  \xi l_i(\xi) \, \diff \xi
    +
    \int_{\xi = -1}^{x_k}  l_i(\xi) \, \diff \xi
  \right]
  .
\end{multline}
Note that the Green matrix is centrosymmetric, hence we need to compute only the half of the elements. Furthermore the first and the last row of the matrix is identically equal to zero. Now our task is to evaluate all integrals in~\eqref{eq:25}. In particular, this requires us to evaluate the integrals of the type
\begin{subequations}
  \label{eq:26}
  \begin{equation}
    \label{eq:27}
    \int_{\xi=-1}^{x_k}  l_i(\xi) \, \diff \xi,
  \end{equation}
  and
  \begin{equation}
    \label{eq:28}
    \int_{\xi=-1}^{x_k}  \xi l_i(\xi) \, \diff \xi,
  \end{equation}
\end{subequations}
for $k, i =0, \dots, N$. Recall that we want exact and fast-to-compute formulae for these integrals,  we do not want to use numerical quadrature or symbolic integration. (Plots of the integrands in~\eqref{eq:26} are shown in Figure~\ref{fig:discrete-dirac}. The integrands are highly oscillatory functions for large values of $N$.) We must propose an approach that works well for large~$N$.

Concerning the second integral~\eqref{eq:28}, we see that the corresponding primitive function can be rewritten as 
\begin{equation}
  \label{eq:29}
  \int  \xi l_i(\xi) \, \diff \xi
  =
  \int  \left(\xi - \xi_i \right) l_i(\xi) \, \diff \xi
  +
  \xi_i
  \int  l_i(\xi) \, \diff \xi
  =
  \bweight_i
  \int  l(\xi) \, \diff \xi
  +
  \xi_i
  \int  l_i(\xi) \, \diff \xi,
\end{equation}
where $\left\{\xi_i \right\}_{i=0}^{N}$ are the collocation points. (We are slightly abusing the notation here---this time the collocation points are denoted by $\xi_i$ instead of $x_k$.) In~\eqref{eq:29} we have exploited the fact that the Lagrange basis $\left\{ l_i \right\}_{i=0}^N$ can be rewritten, see \cite[Chapter 5]{trefethen.ln:approximation}, in the form
\begin{equation}
  \label{eq:30}
  l_j(x)
  =
  l(x)
  \frac{\bweight_j}{x-x_j},
\end{equation}
where
\begin{equation}
  l(x) =_{\bydefinition} \prod_{k=0}^N \left(x-x_k\right)
\end{equation}
denotes the \emph{node polynomial}, and $\left\{ \bweight_j \right\}_{j=0}^N$ denote the \emph{barycentric weights}. The barycentric weights are explicitly known for many polynomial systems. In particular, for the Chebyshev--Gauss--Lobatto points~\eqref{eq:12} we get
\begin{equation}
  \label{eq:31}
  \bweight_i 
  = 
  \begin{cases}
    \frac{1}{2}
    (-1)^i \frac{2^{N-1}}{N}, &i =0, \\ 
    \phantom{\frac{1}{2}}(-1)^i \frac{2^{N-1}}{N}, & i=1, \dots, N-1, \\
    \frac{1}{2}
    (-1)^i \frac{2^{N-1}}{N}, &i =N.
  \end{cases}
\end{equation}
Furthermore, the node polynomial $l(x)$ can be, for the Chebyshev--Gauss--Lobatto points~\eqref{eq:12}, rewritten in terms of Chebyshev polynomials of the first kind,
\begin{equation}
  \label{eq:32}
  l(x) = \frac{1}{2^N} \left( \ChebyshevT{N+1}(x) - \ChebyshevT{N-1}(x) \right),
\end{equation}
see, for example, \cite[Chapter 5]{trefethen.ln:approximation}, and the Chebyshev polynomials are straightforward to integrate,
\begin{equation}
  \label{eq:33}
  \int l(x) \diff x
  =
  \frac{1}{2^{N+1}}
  \left(
    \frac{\ChebyshevT{N+2}(x)}{N+2}
    -
    \frac{2 \ChebyshevT{N}(x)}{N}
    +
    \frac{\ChebyshevT{N-2}(x)}{N-2}
  \right)
\end{equation}
see Section~\ref{sec:chebysh-four-seri} for details. Consequently, we obtain an explicit formula for the first primitive function on the right-hand side of~\eqref{eq:29},
\begin{equation}
  \label{eq:34}
  \int  \xi l_i(\xi) \, \diff \xi
  =
  \bweight_i
  \int  l(\xi) \, \diff \xi
  +
  \xi_i
  \int  l_i(\xi) \, \diff \xi
  =
  \frac{\bweight_i}{2^{N+1}}
  \left(
    \frac{\ChebyshevT{N+2}(x)}{N+2}
    -
    \frac{2 \ChebyshevT{N}(x)}{N}
    +
    \frac{\ChebyshevT{N-2}(x)}{N-2}
  \right)
  +
  \xi_i
  \int  l_i(\xi) \, \diff \xi
  .
\end{equation}
The remaining primitive function on the right-hand side of~\eqref{eq:34} is the same primitive function as the primitive function that appears in integral~\eqref{eq:26}. Consequently our task to evaluate the formula for the Green matrix elements~\eqref{eq:25} \emph{reduces to the evaluation of the integrals}~\eqref{eq:27} and to the \emph{evaluation of Chebyshev polynomials} $\ChebyshevT{N+2}$, $\ChebyshevT{N}$ and $\ChebyshevT{N-2}$ at all collocation points $\left\{ x_k \right\}_{k=0}^N$.

Next we recall that the integrals of the Lagrange basis $\left\{ l_i \right\}_{i=0}^N$ over the whole interval $[-1, 1]$ are known. These are just the quadrature weights $\qweight_i$ for the Clenshaw--Curtis quadrature,
\begin{equation}
  \label{eq:35}
  \qweight_i =_{\bydefinition}
  \int_{\xi=-1}^1
  l_i(\xi)
  \, \diff \xi.
\end{equation}
The quadrature weights can be computed without an explicit calculation of the primitive function and its subsequent evaluation. In principle, the weights are computed via \emph{discrete Fourier transform of a suitably chosen vector}, see~\cite{waldvogel.j:fast}, \cite{trefethen.ln:is}, \cite{sommariva.a:fast}, \cite{trefethen.ln:exactness} and \cite{hale.n.trefethen.ln:chebfun} for further comments on the Clenshaw--Curtis quadrature and the computation of the quadrature weights.
In our case we however do not need only the \emph{integrals over the whole domain} $[-1,1]$, we need more. We need the \emph{integrals over all subintervals} $[x_k, x_l]$ between the corresponding collocation points. Fortunately, as we show below, the integrals $\int_{\xi = -1}^{x_k}  l_i(\xi) \, \diff \xi$ and $\int_{\xi = x_k}^{1}  l_i(\xi) \, \diff \xi$ can also be explicitly evaluated in a numerically feasible way using Chebyshev--Fourier series and the discrete cosine transform.

In principle, it suffices to find the Chebyshev--Fourier series for the Lagrange polynomial $l_i(x)$, integrate in the Chebyshev--Fourier space, transform the result back to the node value representation on a finer set of collocation points, and pick up the right components of this node values vector. The individual steps are described below in Section~\ref{sec:chebysh-four-seri} and Section~\ref{sec:chebysh-four-seri-1}, the final formula is given in Section~\ref{sec:computation--int_xi=}, Algorithm~\ref{alg:discrete-dirac-integrals}.

Once we obtain formulae for the integrals $\int_{\xi = -1}^{x_k}  l_i(\xi) \, \diff \xi$ and $\int_{\xi = x_k}^{1}  l_i(\xi) \, \diff \xi$, we use~\eqref{eq:29} for the integrals $\int_{\xi = -1}^{x_k}  \xi l_i(\xi) \, \diff \xi$ and $\int_{\xi = x_k}^{1}  \xi l_i(\xi) \, \diff \xi$. In evaluating~\eqref{eq:33} at the collocation points we can again use the fact that the values of Chebyshev polynomials at collocation points are fast to find using the discrete Fourier transform, see Section~\ref{sec:eval-chebysh-polyn}, Algorithm~\ref{alg:node-polynomial-int}.

The final step in the Green matrix elements calculation is to collect all the auxilliary results and substitute into~\eqref{eq:25}. This yields
\begin{subequations}
  \label{eq:93}
\begin{equation}
  \label{eq:92}
  \GMatrixc{ki}
  =
  \frac{1}{2}
  \left(
    x_k+1
  \right)
  \left[
    P_{i, N}^{\downarrow}(x_k)
    +
    (x_i-1)
    L_{i, N}^{\downarrow}(x_k)
  \right]
  +
  \frac{1}{2}
  \left(
    x_k-1
  \right)
  \left[
    P_{i, N}^{\uparrow}(x_k)
    +
    (x_i+1)
    L_{i, N}^{\uparrow}(x_k)
  \right],
\end{equation}
where we have denoted
\begin{align}
  \label{eq:94}
  P_{i, N}^{\downarrow}(\zeta) &=_{\bydefinition}
                            \left.
                            \frac{\bweight_i}{2^{N+1}}
                            \left(
                            \frac{\ChebyshevT{N+2}(x)}{N+2}
                            -
                            \frac{2 \ChebyshevT{N}(x)}{N}
                            +
                            \frac{\ChebyshevT{N-2}(x)}{N-2}
                            \right)
                            \right|_{x = 1}
                            -
                            \left.
                            \frac{\bweight_i}{2^{N+1}}
                            \left(
                            \frac{\ChebyshevT{N+2}(x)}{N+2}
                            -
                            \frac{2 \ChebyshevT{N}(x)}{N}
                            +
                            \frac{\ChebyshevT{N-2}(x)}{N-2}
                            \right)
                            \right|_{x = \zeta}
                            ,
  \\
  \label{eq:95}
  P_{i, N}^{\uparrow}(\zeta) &=_{\bydefinition}
                           \left.
                            \frac{\bweight_i}{2^{N+1}}
                            \left(
                            \frac{\ChebyshevT{N+2}(x)}{N+2}
                            -
                            \frac{2 \ChebyshevT{N}(x)}{N}
                            +
                            \frac{\ChebyshevT{N-2}(x)}{N-2}
                            \right)
                            \right|_{x = \zeta}
                            -
                            \left.
                            \frac{\bweight_i}{2^{N+1}}
                            \left(
                            \frac{\ChebyshevT{N+2}(x)}{N+2}
                            -
                            \frac{2 \ChebyshevT{N}(x)}{N}
                            +
                            \frac{\ChebyshevT{N-2}(x)}{N-2}
                            \right)
                          \right|_{x = -1}
                          ,
  \\
  \label{eq:96}
  L_{i, N}^{\downarrow}(\zeta) & =_{\bydefinition}\int_{\xi = \zeta}^{1}  l_i(\xi) \, \diff \xi, \\
  \label{eq:97}
  L_{i, N}^{\uparrow}(\zeta) &=_{\bydefinition} \int_{\xi = -1}^{\zeta}  l_i(\xi) \, \diff \xi.
\end{align}
\end{subequations}
Formula~\eqref{eq:92} is the final formula in Algorithm~\ref{alg:green-matrix-elements} for fast construction of Green matrix.

Before we proceed with the detailed description of the algorithms, we briefly recall basic properties of Chebyshev polynomials, in particular the identities for primitive functions, see Section~\eqref{sec:chebysh-four-seri}, the Chebyshev--Fourier series and node values/Fourier coefficients conversion via the discrete cosine transformation, see Section~\eqref{sec:chebysh-four-seri-1}. 

\subsection{Chebyshev--Fourier series and primitive functions}
\label{sec:chebysh-four-seri}
Let us consider a function $u$ that is a polynomial of degree at most $N$. The Chebyshev--Fourier series of $u$ is the series 
\begin{equation}
  \label{eq:38}
   u(x) = \sum_{j=0}^{N} \widehat{u}_j \ChebyshevT{j}(x),
 \end{equation}
 where $\ChebyshevT{j}$ denote the Chebyshev polynomials of the first kind. We note that on the discrete level we record only the coefficients $\left\{ \widehat{u}_j \right\}_{j=0}^N$, i.e. we work with vectors
 \begin{equation}
   \label{eq:39}
   \widehat{\vec{u}}
   =
   \transpose{
   \begin{bmatrix}
     \widehat{u}_0 &
     \widehat{u}_1 &
     \widehat{u}_2 &
     \cdots
     &
     \widehat{u}_j
     &
     \cdots
     &
     \widehat{u}_{N-1} &
     \widehat{u}_{N}
   \end{bmatrix}
   }
   ,
 \end{equation}
 and we recall that for $u$ being a polynomial of degree at most $N$ the corresponding Chebyshev--Fourier representation via the coefficients  $\left\{ \widehat{u}_j \right\}_{j=0}^N$ is exact.

 If we know the Chebyshev--Fourier series for $u$, then it is straightforward to find the primitive function to $u$.   The Chebyshev polynomials are integrated according to the formulae
 \begin{subequations}
   \label{eq:40}
   \begin{align}
     \label{eq:41}
     \int \ChebyshevT{0}(\xi) \, \diff \xi &= x, \\
     \label{eq:42}
     \int \ChebyshevT{1}(\xi) \, \diff \xi &= \frac{x^2}{2}, \\
     \label{eq:43}
     \int
     \ChebyshevT{j}(\xi)
     \, \diff \xi
     &=
     \frac{1}{2}
     \left(
     \frac{\ChebyshevT{j+1}(\xi)}{j+1}
     -
     \frac{\ChebyshevT{j-1}(\xi)}{j-1}
       \right)
       , \qquad j \in \N,\ j \geq 2,
   \end{align}
 \end{subequations}
 which can be by virtue of relations $x = \ChebyshevT{1}(x)$ and $\frac{x^2}{2} = \frac{\ChebyshevT{2}(x) + \ChebyshevT{0}(x)}{4}$ further rewritten as
 \begin{subequations}
   \label{eq:44}
   \begin{align}
     \label{eq:45}
     \int \ChebyshevT{0}(\xi) \, \diff \xi &= \ChebyshevT{1}(x), \\
     \label{eq:46}
     \int \ChebyshevT{1}(\xi) \, \diff \xi &= \frac{\ChebyshevT{2}(x) + \ChebyshevT{0}(x)}{4}, \\
     \label{eq:47}
     \int
     \ChebyshevT{j}(\xi)
     \, \diff \xi
     &=
     \frac{1}{2}
     \left(
     \frac{\ChebyshevT{j+1}(\xi)}{j+1}
     -
     \frac{\ChebyshevT{j-1}(\xi)}{j-1}
       \right)
       , \qquad j \in \N,\ j \geq 2,
   \end{align}
 \end{subequations}
 see, for example, \cite{mason.jc.handscomb.dc:chebyshev} or \cite{trefethen.ln:approximation}.
 
 Consequently, the Chebyshev--Fourier series for the primitive function $\int u(\xi) \, \diff \xi$ is easy to write down. (The primitive function of course includes some integration constant that is fixed later.) \emph{The polynomial degree of the primitive function to~$u$ is however one degree higher than the polynomial degree of the original function $u$}. Consequently, in the discrete setting we must work with longer vectors/larger number of coefficients in the series~\eqref{eq:38}.

 It turns out that it is convenient to go from vectors of dimension $N+1$, that is from $N+1$ collocation points and $N$ degree polynomials, to vectors of dimension $2N+1$ and polynomials of degree $2N$. This allows us to find the primitive function to \eqref{eq:38} as follows. First we take the vector $\widehat{\vec{u}}$ of length $N+1$ and we extend it by $N$ trailing zeros. This generates vector of length $2N+1$, which is then integrated using rules~\eqref{eq:44}. This yields
 \begin{equation}
   \label{eq:48}
   \begin{bmatrix}
     \widehat{u}_0 \\
     \widehat{u}_1 \\
     \widehat{u}_2 \\
     \vdots
     \\
     \widehat{u}_j
     \\
     \vdots
     \\
     \widehat{u}_{N-1} \\
     \widehat{u}_{N}
   \end{bmatrix}
   \stackrel{\text{extend}}{\mapsto}
   \begin{bmatrix}
     \widehat{u}_0 \\
     \widehat{u}_1 \\
     \widehat{u}_2 \\
     \vdots
     \\
     \widehat{u}_j
     \\
     \vdots
     \\
     \widehat{u}_{N-1} \\
     \widehat{u}_{N} \\
     0 \\
     0 \\
     \vdots \\
     0
   \end{bmatrix}
   \stackrel{\text{integrate}}{\mapsto}
   \begin{bmatrix}
     \frac{\widehat{u}_1}{4} \\
     \widehat{u}_0 - \frac{\widehat{u}_2}{2} \\
     \frac{1}{4} \left( \widehat{u}_{1} - \widehat{u}_{3} \right) \\
     \vdots \\
     \frac{1}{2j} \left( \widehat{u}_{j-1} - \widehat{u}_{j+1} \right) \\
     \vdots \\
     \frac{1}{2(N-1)} \left( \widehat{u}_{N-2} - \widehat{u}_{N} \right) \\
     \frac{1}{2N} \left( \widehat{u}_{N} - \widehat{u}_{N+1} \right) \\
     \frac{1}{2(N+1)} \left( \widehat{u}_{N} - \widehat{u}_{N+2} \right) \\
     0 \\
     \vdots \\
     0
   \end{bmatrix}
   .
 \end{equation}
 (The original vector $\widehat{\vec{u}}$ has dimension $N+1$, and its last element is $\widehat{u}_N$. The elements $\widehat{u}_{N+1}$ and $\widehat{u}_{N+2}$ in fact belong to the extended vector, and they are conveniently set equal to zero. This later enables us to write the integration operation using one formula.) The first two elements of the final vector are exceptional in the sense that they must be computed using a special rule. The remaining elements can be calculated using the simple rule $\frac{1}{2j} \left( \widehat{u}_{j-1} - \widehat{u}_{j+1} \right)$, which gives the coefficient for the Chebyshev polynomial $\ChebyshevT{j}(x)$, that is for the $(j-1)$-th row of the corresponding vector provided that we number the elements from one. Note that~\eqref{eq:48} in fact gives us the primitive function $\int_{\xi=x_N}^x l_i(\xi) \, \diff \xi + C$, where $C$ is the integration constant; the lower bound of the integral is, without loss of generality, conveniently chosen as $x_N=-1$. 

 For later reference we denote the extension operator as $\EXTEND$ and the discrete integration operator as $\INTEGRATE$, that is we define
 \begin{subequations}
   \label{eq:49}
   \begin{equation}
     \label{eq:50}
     \EXTEND \left(\widehat{\vec{u}}\right)
     =_{\bydefinition}
       \begin{bmatrix}
         \widehat{\vec{u}} \\
         \vec{0}_N
       \end{bmatrix}
     ,
   \end{equation}
   and
   \begin{equation}
   \label{eq:51}
   \INTEGRATE
   \left(
     \begin{bmatrix}
       \widehat{\vec{u}} \\
       \vec{0}_N
     \end{bmatrix}
   \right)
   =_{\bydefinition}
     \begingroup
     \renewcommand*{\arraystretch}{1.3}
     \begin{bmatrix}
     \frac{\widehat{u}_1}{4} \\
     \widehat{u}_0 - \frac{\widehat{u}_2}{2} \\
     \frac{1}{4} \left( \widehat{u}_{1} - \widehat{u}_{3} \right) \\
     \frac{1}{6} \left( \widehat{u}_{2} - \widehat{u}_{4} \right) \\
     \frac{1}{8} \left( \widehat{u}_{3} - \widehat{u}_{5} \right) \\
     \vdots \\
     \frac{1}{2j} \left( \widehat{u}_{j-1} - \widehat{u}_{j+1} \right) \\
     \vdots \\
     \frac{1}{2(N-1)} \left( \widehat{u}_{N-2} - \widehat{u}_{N} \right) \\
     \frac{1}{2N} \left( \widehat{u}_{N-1} - \widehat{u}_{N+1} \right) \\
     \frac{1}{2(N+1)} \left( \widehat{u}_{N} - \widehat{u}_{N+2} \right) \\
     0 \\
     0 \\
     \vdots \\
     0
   \end{bmatrix}
   \endgroup
   =
     \begingroup
     \renewcommand*{\arraystretch}{1.3}
     \tensorschur{
     \begin{bmatrix}
       0 \\
       \frac{1}{2} \\
       \frac{1}{4} \\
       \frac{1}{6} \\
       \frac{1}{8} \\
       \vdots \\
       \frac{1}{2j} \\
       \vdots \\
       \frac{1}{2(N-1)} \\
       \frac{1}{2N} \\
       \frac{1}{2(N+1)} \\
       0 \\
       0 \\
       \vdots \\
       0
     \end{bmatrix}
     }
     {
     \left(
     \begin{bmatrix}
       0 \\
       \widehat{u}_0 \\
       \widehat{u}_1 \\
       \widehat{u}_2 \\
       \widehat{u}_3 \\
       \vdots \\
       \widehat{u}_{j-1} \\
       \vdots \\
       \widehat{u}_{N-2}\\
       \widehat{u}_{N-1} \\
       \widehat{u}_{N} \\
       0 \\
       0 \\
       \vdots \\
       0
     \end{bmatrix}
     -
     \begin{bmatrix}
       \widehat{u}_1 \\
       \widehat{u}_2 \\
       \widehat{u}_3 \\
       \widehat{u}_4 \\
       \widehat{u}_5 \\
       \vdots \\
       \widehat{u}_{j+1} \\
       \vdots \\
       \widehat{u}_{N} \\
       0 \\
       0 \\
       0 \\
       0 \\
       \vdots \\
       0
     \end{bmatrix}
     \right)
     }
     +
     \begin{bmatrix}
       \frac{\widehat{u}_1}{4} \\
       \frac{\widehat{u_0}}{2} \\
       0 \\
       0 \\
       0 \\
       \vdots \\
       0 \\
       \vdots \\
       0 \\
       0 \\
       0 \\
       0 \\
       0 \\
       \vdots \\
       0
     \end{bmatrix}
     \endgroup
     ,
 \end{equation}%
\end{subequations}
where $\vec{0}_N$ denotes the zero vector of dimension $N$, and where the multiplication $\tensorschur{}{}$ on the right-hand side is the element-wise product (Schur product, Hadamard product). The two vectors on the right-hand side of~\eqref{eq:51} can be obtained from the (extended) original vector $\widehat{\vec{u}}$, see~\eqref{eq:39}, by the standard circular shift operations, see~\texttt{circshift} in Matlab and \texttt{RotateLeft}/\texttt{RotateRight} in Wolfram Mathematica. Note that this operator represents the integration provided that the input has a particular structure with trailing zeros.

\subsection{Chebyshev--Fourier series, node values, and the discrete cosine transform}
\label{sec:chebysh-four-seri-1}
For function $u$ being a polynomial of degree at most $N$ there is one-to-one relation between the node values $\left\{ u_i \right\}_{i=0}^N$ at the Chebyshev points and the Chebyshev--Fourier coefficients $\left\{\widehat{u}_i\right\}_{i=0}^{N}$. We can \emph{equivalently} use either the node value representation
\begin{equation}
  \label{eq:52}
   u(x) = \sum_{i=0}^N u_i l_i(x)
 \end{equation}
 or the Chebyshev--Fourier representation
 \begin{equation}
   \label{eq:53}
   u(x) = \sum_{j=0}^{N} \widehat{u}_j \ChebyshevT{j}(x).
 \end{equation}
 On the discrete level we again represent~\eqref{eq:52} using the vector of node values
 \begin{equation}
   \label{eq:54}
   \vec{u}
   =
   \transpose{
   \begin{bmatrix}
     u_0 &
     u_1 &
     u_2 &
     \cdots
     &
     u_j
     &
     \cdots
     &
     u_{N-1} &
     u_{N}
   \end{bmatrix}
 }
 \end{equation}

 Since the Chebyshev polynomials of the first kind are defined as
 $
 \ChebyshevT{k}(x) = _{\bydefinition} \cos \left( k \arccos x \right)
 $, we see that if we evaluate the $k$-th Chebyshev polynomial $\ChebyshevT{k}$ at the $i$-th collocation point $x_i$, see~\eqref{eq:12}, we get
 \begin{equation}
   \label{eq:55}
   \ChebyshevT{k}(x_i)
   =
   \cos \left( k \frac{i}{N} \pi \right).
 \end{equation}
 Consequently, if we use~\eqref{eq:55} in~\eqref{eq:53}, we see that
 \begin{equation}
   \label{eq:56}
   u_i
   =
   \sum_{j=0}^{N} \widehat{u}_j \cos \left( j \frac{i}{N} \pi \right)
   =
   \widehat{u}_0
   +
   \sum_{j=1}^{N-1} \widehat{u}_j \cos \left( j \frac{i}{N} \pi \right)
   +
   (-1)^i
   \widehat{u}_{N}
   .
 \end{equation}
 This expression can be compared with the formula for the discrete cosine transform, and it turns out that the transformation from the node value representation $\vec{u}$ and the Chebyshev--Fourier representation $\widehat{\vec{u}}$ and vice versa can be done via the discrete cosine transform, see, for example, \cite[Chapter 8]{trefethen.ln:spectral}, \cite[Section 4.7]{mason.jc.handscomb.dc:chebyshev} and~\cite[Section 3.4.1]{shen.j.tang.t.ea:spectral}. Thus we can switch between the representations in a numerically convenient---fast and stable---way.

 The definitions of discrete cosine transform differ in various software packages; different normalisations and scaling factors are used. To be specific, we describe definitions of discrete cosine transform in Matlab and Wolfram Mathematica.

\subsubsection{Discrete cosine transform---Wolfram Mathematica}
\label{sec:discr-cosine-transf-1}
 The DCT-I variant of the discrete cosine transform of vector $\vec{v}$ of dimension $n$ with elements numbered by $s=1, \dots, n$ is defined as
 \begin{equation}
   \label{eq:57}
   \widehat{v}_\vecelem{s}
   =
   \sqrt{\frac{2}{n-1}}
   \left(
     \frac{v_{\vecelem{1}}}{2}
     +
     \sum_{r=2}^{n-1} v_{\vecelem{r}} \cos \left( \frac{\pi}{n-1}  \left(r-1 \right) \left(s-1\right)  \right)
     +
     \left(-1 \right)^{s-1} \frac{v_{\vecelem{n}}}{2}
   \right),
 \end{equation}
 where we have used the convention implemented in Wolfram Mathematica 
in the function \texttt{FourierDCT}. If we use the definition~\eqref{eq:57}, then the inverse transformation to DCT-I is again DCT-I. Inspecting the similarity between~\eqref{eq:56} and~\eqref{eq:57}, we see that we can write~\eqref{eq:56} as
 \begin{equation}
   \label{eq:58}
   u_{i-1}
   =
   \frac{2 \widehat{u}_0}{2}
   +
   \sum_{j=2}^{N} \widehat{u}_{j-1} \cos \left( \frac{\pi}{N}  (j-1) (i-1)\right)
   +
   (-1)^{i-1}
   \frac{2 \widehat{u}_{N}}{2}
 \end{equation}
 for $i=1, \dots, N+1$, which means that the vector of node values $\vec{u}$ is obtained by the discrete cosine transform of the suitably modified vector of Chebyshev--Fourier coefficients $\widehat{\vec{u}}$,  
 \begin{equation}
   \label{eq:59}
   \begin{bmatrix}
     u_0 \\
     u_1 \\
     u_2 \\
     \vdots
     \\
     u_j
     \\
     \vdots
     \\
     u_{N-1} \\
     u_{N}
   \end{bmatrix}
   =
   \DCTI_{\text{Wolfram}}
   \left(
     \sqrt{\frac{N}{2}}
     \tensorschur{
       \begin{bmatrix}
         2 \\
         1 \\
         1 \\
         \vdots \\
         1 \\
         \vdots \\
         1 \\
         2
       \end{bmatrix}
     }
     {
       \begin{bmatrix}
         \widehat{u}_0 \\
         \widehat{u}_1 \\
         \widehat{u}_2 \\
         \vdots
         \\
         \widehat{u}_j
         \\
         \vdots
         \\
         \widehat{u}_{N-1} \\
         \widehat{u}_{N}
       \end{bmatrix}
     }
   \right)
   .
 \end{equation}
 Transformation~\eqref{eq:59} from Chebyshev--Fourier coefficients to node values can be rewritten as
 \begin{equation}
   \label{eq:60}
   \vec{u}
   =
   \DCTI_{\text{Wolfram}}
   \left(
     \sqrt{\frac{N}{2}}
     \tensorschur{
       \vec{c}
     }
     {
       \widehat{\vec{u}}
     }
   \right),
 \end{equation}
 while the inverse transformation from node values to Chebyshev--Fourier coefficients reads
 \begin{equation}
   \label{eq:61}
   \widehat{\vec{u}}
   =
   \tensorschur{
     \vec{c}^{-1}
   }
   {
     \DCTI_{\text{Wolfram}}
   \left(
     \sqrt{\frac{2}{N}}
     \vec{u}
   \right)
 },
\end{equation}
where we use notation $\vec{c}$ for the vector of normalisation coefficients,
\begin{equation}
  \label{eq:62}
  \vec{c}
  =_{\bydefinition}
  \transpose{
  \begin{bmatrix}
    2 &
    1 &
    1 &
    \cdots &
    1 &
    1 &
    2
  \end{bmatrix}
  }
  ,
  \qquad
  \vec{c}^{-1}
  =_{\bydefinition}
  \transpose{
    \begin{bmatrix}
    \frac{1}{2} &
    1 &
    1 &
    \cdots &
    1 &
    1 &
    \frac{1}{2}
  \end{bmatrix}
  }
  .
\end{equation}
Note that we use the same notation for the vector with different dimensions, that is $\vec{c}$ can have the dimension $N+1$ or $2N+1$ depending on the context. 

\subsubsection{Discrete cosine transform---Matlab}
\label{sec:discr-cosine-transf-2}
 The DCT-I variant of the discrete cosine transform of vector $\vec{v}$ of dimension $n$ with elements numbered by $s=1, \dots, n$ is defined as
 \begin{equation}
   \label{eq:63}
   \widehat{v}_\vecelem{s}
   =
   \begin{cases}
     \frac{1}{\sqrt{2}}
     \sqrt{\frac{2}{n-1}}
     \left(
       \frac{v_{\vecelem{1}}}{\sqrt{2}}
       +
       \sum_{r=2}^{n-1} v_{\vecelem{r}} \cos \left( \frac{\pi}{n-1}  \left(r-1 \right) \left(s-1\right)  \right)
       +
       \left(-1 \right)^{s-1} \frac{v_{\vecelem{n}}}{\sqrt{2}}
     \right),
     &
     s = 1, 
     \\
     \phantom{\frac{1}{\sqrt{2}}}
     \sqrt{\frac{2}{n-1}}
     \left(
       \frac{v_{\vecelem{1}}}{\sqrt{2}}
       +
       \sum_{r=2}^{n-1} v_{\vecelem{r}} \cos \left( \frac{\pi}{n-1}  \left(r-1 \right) \left(s-1\right)  \right)
       +
       \left(-1 \right)^{s-1} \frac{v_{\vecelem{n}}}{\sqrt{2}}
     \right),
     &
     s=2, \dots, n-1, \\
     \frac{1}{\sqrt{2}}
     \sqrt{\frac{2}{n-1}}
     \left(
       \frac{v_{\vecelem{1}}}{\sqrt{2}}
       +
       \sum_{r=2}^{n-1} v_{\vecelem{r}} \cos \left( \frac{\pi}{n-1}  \left(r-1 \right) \left(s-1\right)  \right)
       +
       \left(-1 \right)^{s-1} \frac{v_{\vecelem{n}}}{\sqrt{2}}
     \right),
     & s = n,
   \end{cases}
 \end{equation}
 where we have used the convention implemented in Matlab in the function \texttt{dct}. If we use the definition~\eqref{eq:63}, then the inverse transformation to DCT-I is again DCT-I. Inspecting the similarity between~\eqref{eq:56} and~\eqref{eq:63}, we see that we can write~\eqref{eq:56} as
 \begin{equation}
   \label{eq:64}
   u_{i-1}
   =
   \frac{\sqrt{2} \widehat{u}_0}{\sqrt{2}}
   +
   \sum_{j=2}^{N} \widehat{u}_{j-1} \cos \left( \frac{\pi}{N}  (j-1) (i-1)\right)
   +
   (-1)^{i-1}
   \frac{\sqrt{2} \widehat{u}_{N}}{\sqrt{2}}
   ,
 \end{equation}
 for $i=1, \dots, N+1$, which means that the vector of node values $\vec{u}$ is obtained by the discrete cosine transform of the suitably modified vector of Chebyshev--Fourier coefficients $\widehat{\vec{u}}$,  
 \begin{equation}
   \label{eq:157a}
   \begin{bmatrix}
     u_0 \\
     u_1 \\
     u_2 \\
     \vdots
     \\
     u_j
     \\
     \vdots
     \\
     u_{N-1} \\
     u_{N}
   \end{bmatrix}
   =
   \tensorschur{
     \begin{bmatrix}
       \sqrt{2} \\
       1 \\
       1 \\
       \vdots \\
       1 \\
       \vdots \\
       1 \\
       \sqrt{2}
     \end{bmatrix}
   }
   {
   \DCTI_{\text{Matlab}}
   \left(
     \sqrt{\frac{N}{2}}
     \tensorschur{
       \begin{bmatrix}
         \sqrt{2} \\
         1 \\
         1 \\
         \vdots \\
         1 \\
         \vdots \\
         1 \\
         \sqrt{2}
       \end{bmatrix}
     }
     {
       \begin{bmatrix}
         \widehat{u}_0 \\
         \widehat{u}_1 \\
         \widehat{u}_2 \\
         \vdots
         \\
         \widehat{u}_j
         \\
         \vdots
         \\
         \widehat{u}_{N-1} \\
         \widehat{u}_{N}
       \end{bmatrix}
     }
   \right)
   }
 \end{equation}
 where the multiplication $\tensorschur{}{}$ on the right-hand side is again the element-wise product (Schur product, Hadamard product). Transformation~\eqref{eq:157a} from Chebyshev--Fourier coefficients to node values can be rewritten as
 \begin{equation}
   \label{eq:158a}
   \vec{u}
   =
   \tensorschur{
     \sqrt{
         \vec{c}
       }
     }
     {
   \DCTI_{\text{Matlab}}
   \left(
     \sqrt{\frac{N}{2}}
     \tensorschur{
       \sqrt{
         \vec{c}
       }
     }
     {
       \widehat{\vec{u}}
     }
   \right)
   },
 \end{equation}
 while the inverse transformation from node values to Chebyshev--Fourier coefficients reads
 \begin{equation}
   \label{eq:97a}
   \widehat{\vec{u}}
   =
   \tensorschur{
     \left(\sqrt{\vec{c}}\right)^{-1}
   }
   {
     \DCTI_{\text{Matlab}}
   \left(
     \sqrt{\frac{2}{N}}
     \left(\sqrt{\vec{c}}\right)^{-1}
     \vec{u}
   \right)
 },
\end{equation}
where we use notation $\vec{c}$ for the vectors of normalisation coefficients, see~\eqref{eq:62}, and where the square root of this vector is understood as the application of the square root to every element of the vector.

\subsection{Computation of $\int_{\xi=-1 }^{x_k} l_i(\xi) \, \diff \xi$ and $\int_{\xi=x_k }^{1} l_i(\xi) \, \diff \xi$}
\label{sec:computation--int_xi=}
The Lagrange basis polynomial $l_i$ is easy to write down using the node value representation. Indeed, the corresponding vector is
\begin{equation}
  \label{eq:69}
  \vec{l}_i
  =
  \transpose{
  \begin{bmatrix}
    0 &
    \cdots &
    0 &
    1 &
    0 &
    \cdots &
    0
  \end{bmatrix}
  }
  ,
\end{equation}
where the only non-zero element is on the $(i+1)$-th position, provided that we number the vector elements from one. This is a vector of dimension $N+1$. Now we want to find the formula for integrals
\begin{equation}
  \label{eq:70}
  \int_{\xi=x_l }^{x_k} l_i(\xi) \, \diff \xi,
\end{equation}
where $\left\{ x_i \right\}_{i=0}^N$ denote the Chebyshev--Gauss--Lobatto points. In order to do that, we first transform~\eqref{eq:69} to the coefficient space, see Section~\ref{sec:chebysh-four-seri-1},
\begin{subequations}
  \label{eq:71}
  \begin{align}
    \label{eq:72}
    \widehat{\vec{l}}_i
    &=
      \tensorschur{
      \vec{c}^{-1}
      }
      {
      \DCTI_{\text{Wolfram}}
      \left(
      \sqrt{\frac{2}{N}}
      \vec{l}_i
      \right)
      }
      ,
    \\
    \label{eq:73}
    \widehat{\vec{l}}_i
    &=
    \tensorschur{
    \left(\sqrt{\vec{c}}\right)^{-1}
    }
    {
    \DCTI_{\text{Matlab}}
    \left(
    \sqrt{\frac{2}{N}}
    \left(\sqrt{\vec{c}}\right)^{-1}
    \vec{l}_i
    \right)
    }
    ,
  \end{align}
\end{subequations}
where we use the formula that corresponds to our choice regarding the definition of discrete cosine transform. In this case the discrete cosine transform operates on the vector of dimension~$N+1$ and $ \vec{c}^{-1}$ is understood as the vector of dimension~$N+1$.

Once we are done with the transformation to the Chebyshev--Fourier space, we extend the space for Chebyshev--Fourier coefficients, in particular we add $N$ trailing zeros $ \vec{0}_{N}$ to the column vector $\widehat{\vec{l}}_i $. (Recall that the integration increases the polynomial order by one, hence we need more Chebyshev polynomials in the Chebyshev--Fourier expansion.) This can be encoded as an application of the operator $\EXTEND$,
\begin{equation}
  \label{eq:74}
  \EXTEND \left( \widehat{\vec{l}}_i \right)
  =
  \begin{bmatrix}
    \widehat{\vec{l}}_i \\
    \vec{0}_{N}
  \end{bmatrix}.
\end{equation}
We denote the outcome of this operation as $\widehat{\vec{l}}_{i, 2N}$, that is we set
\begin{equation}
  \label{eq:98}
  \widehat{\vec{l}}_{i, 2N} =_{\bydefinition} \EXTEND \left( \widehat{\vec{l}}_i \right),
\end{equation}
and we recall that $\widehat{\vec{l}}_{i, 2N}$ is a vector of dimension~$2N+1$.

Now we are in position to find the function
\begin{equation}
  \label{eq:75}
  L_{i, N}(x) =_{\bydefinition} \int_{\xi=x_N}^x l_i(\xi) \, \diff \xi + C,
\end{equation}
which we shortly refer to as the primitive function. (The index $N$ refers to the fact that we are integrating Lagrange basis for the degree~$N$ polynomial interpolation. It does not refer to the dimension of a vector,  $L_{i, N}$ denotes a function.) The Chebyshev--Fourier representation of the primitive function is given by the formula
\begin{equation}
  \label{eq:76}
  \INTEGRATE
  \left(
    \widehat{\vec{l}}_{i, 2N}
  \right)
  .
\end{equation}
\emph{Regarding the computation of Green matrix elements, we are, however, not interested in the primitive function itself, we want to \emph{evaluate} the primitive function at collocation points $\left\{x_k\right\}_{k=0}^N$}. This is easy to do, if we transform from the Chebyshev--Fourier representation back to the node value representation.

The transformation of
$
\INTEGRATE
\left(
  \widehat{\vec{l}}_{i, 2N}
\right)
$
back to the node value representation is done using the discrete cosine transform,
\begin{subequations}
  \label{eq:77}
  \begin{align}
    \label{eq:78}
  \vec{L}_{i, 2N}
  &=
  \DCTI_{\text{Wolfram}}
   \left(
     \sqrt{N}
     \tensorschur{
       \vec{c}
     }
     {
       \INTEGRATE
       \left(
         \widehat{\vec{l}}_{i, 2N}
       \right)
     }
    \right), \\
    \label{eq:79}
  \vec{L}_{i, 2N}
  &=
    \tensorschur{
     \sqrt{
         \vec{c}
       }
     }
     {
   \DCTI_{\text{Matlab}}
   \left(
     \sqrt{N}
     \tensorschur{
       \sqrt{
         \vec{c}
       }
     }
     {
       \INTEGRATE
       \left(
         \widehat{\vec{l}}_{i, 2N}
       \right)
     }
   \right)
    }
    ,
  \end{align}
\end{subequations}
and we obtain the vector $\vec{L}_{i, 2N}$ of values of the primitive function $ L_{i, N}$ at the Chebyshev--Gauss--Lobatto points~$\left\{ y_j \right\}_{j=0}^{2N}$,
 \begin{equation}
   \label{eq:80}
   y_k = _{\bydefinition} \cos \left( \frac{k \pi}{2N}  \right).
 \end{equation}
 (We note that the discrete cosine transform in~\eqref{eq:77} is now applied on the vector of dimension $2N+1$, and that $\vec{c}$ is understood as the vector of dimension $2N+1$.  The choice of formula in~\eqref{eq:77} again depends on the particular definition of the discrete cosine transform we are using.) This gives us the values of $L_{i, N}$, \emph{but on a finer grid}.

 In other words, we have at our disposal the values
 \begin{equation}
   \label{eq:81}
   \vec{L}_{i, 2N}
   =
   \begin{bmatrix}
     L_{i, N}(y_0) \\
     L_{i, N}(y_1) \\
     L_{i, N}(y_2) \\
     \vdots \\
     L_{i, N}(y_{2N-1}) \\
     L_{i, N}(y_{2N}) \\
   \end{bmatrix}
   ,
 \end{equation}
but we want the values at the original collocation points $\left\{ x_k \right\}_{k=0}^N$. However, due to the convenient choice of the finer grid, we have the interlacing property between the original  Chebyshev--Gauss--Lobatto points~$\left\{ x_k \right\}_{k=0}^N$  and the new  Chebyshev--Gauss--Lobatto points~$\left\{ y_j \right\}_{j=0}^{2N}$. In particular, for $m=0, \dots, N$ we have
 \begin{equation}
   \label{eq:82}
   x_m = y_{2m},
 \end{equation}
 which implies
 \begin{equation}
   \label{eq:83}
   L_i(x_m)
   =
   L_i(y_{2m})
 \end{equation}
 Consequently, if we take every other element from the vector $\vec{L}_{i, 2N}$, then we get the vector of values at the original collocation points. This can be written as
 \begin{equation}
   \label{eq:84}
   \vec{L}_{i, N}
   =_{\bydefinition}
   \REDUCE
   \left(
     \vec{L}_{i, 2N}
   \right),
 \end{equation}
 where the reduction operator $\REDUCE$ is defined as
 \begin{equation}
   \label{eq:85}
   \REDUCE
   \left(
     \begin{bmatrix}
       a_0 \\
       a_1 \\
       a_2 \\
       a_3 \\
       a_4 \\
       \vdots \\
       a_{2N-2} \\
       a_{2N-1} \\
       a_{2N}
     \end{bmatrix}
   \right)
   =_{\bydefinition}
     \begin{bmatrix}
       a_0 \\
       a_2 \\
       a_4 \\
       \vdots \\
       a_{2N-2} \\
       a_{2N}
     \end{bmatrix}
     .
   \end{equation}
   Note that at this point we lose the one-to-one correspondence between the node values and the Chebyshev--Fourier coefficients. The node values at the original set of collocation points are \emph{exact} as desired, but they are, due to the integration, node values for polynomial of degree $N+1$, hence the corresponding Chebyshev--Fourier series constructed from the node values on the original set of $N+1$ collocation points would be a \emph{truncated} Chebyshev--Fourier series only. But since we need only the node values, this is not a problem. 
   
 The integrals
$
   \int_{\xi=x_l }^{x_k} l_i(\xi) \, \diff \xi
$
 are then computed by subtracting the corresponding elements in the reduced vector $\vec{L}_{i, N}$.  Since we are interested in the integrals $\int_{\xi=-1 }^{x_k} l_i(\xi) \, \diff \xi = \int_{\xi=x_N }^{x_k} l_i(\xi) \, \diff \xi$, we see that the last operation we need is
 \begin{equation}
   \label{eq:86}
   \vec{L}_{i, N}^{\uparrow}
   =
   \vec{L}_{i, N}
   -
   L_{i, N}(x_N)
   \vec{1}_{N+1},
 \end{equation}
 where $L_i(x_N)$ is the last element in vector $\vec{L}_{i, N}$ and $\vec{1}_{N+1}$ denotes the vector of ones with dimension $N+1$. This operation fixes the integration constant to the right value. Finally, we get the vector
 \begin{equation}
   \label{eq:87}
   \vec{L}_{i, N}^{\uparrow}
   =
   \begin{bmatrix}
     L_{i, N}^{\uparrow}(x_0) \\
     L_{i, N}^{\uparrow}(x_1) \\
     L_{i, N}^{\uparrow}(x_2) \\
     \vdots \\
     L_{i, N}^{\uparrow}(x_{N-1}) \\
     L_{i, N}^{\uparrow}(x_{N}) \\
   \end{bmatrix}
   =
   \begin{bmatrix}
     \int_{\xi=-1 }^{x_0} l_i(\xi) \, \diff \xi \\
     \int_{\xi=-1 }^{x_1} l_i(\xi) \, \diff \xi \\
     \int_{\xi=-1 }^{x_2} l_i(\xi) \, \diff \xi \\
     \vdots \\
     \int_{\xi=-1 }^{x_{N-1}} l_i(\xi) \, \diff \xi \\
     \int_{\xi=-1 }^{x_{N}} l_i(\xi) \, \diff \xi
   \end{bmatrix}
   .
 \end{equation}
 The last element $\int_{\xi=-1 }^{x_{N}} l_i(\xi) \, \diff \xi$ must be equal to zero, which is granted by the proposed construction/choice of the integration constant.

 The integrals $\int_{\xi=x_k }^{1} l_i(\xi) \, \diff \xi = \int_{\xi=x_k }^{x_0} l_i(\xi) \, \diff \xi$ are computed in a similar manner, it suffices to appropriately fix the integration constant  
 \begin{equation}
   \label{eq:88}
   \vec{L}_{i, N}^{\downarrow}
   =
   L_{i, N}(x_0)
   \vec{1}_{N+1}
   -
   \vec{L}_{i, N}
   ,
 \end{equation}
 and the vector $\vec{L}_{i, N}^{\downarrow}$ contains the values
 \begin{equation}
   \label{eq:89}
   \vec{L}_{i, N}^{\downarrow}
   =
     \begin{bmatrix}
     L_{i, N}^{\downarrow}(x_0) \\
     L_{i, N}^{\downarrow}(x_1) \\
     L_{i, N}^{\downarrow}(x_2) \\
     \vdots \\
     L_{i, N}^{\downarrow}(x_{N-1}) \\
     L_{i, N}^{\downarrow}(x_{N}) \\
   \end{bmatrix}
   =
   \begin{bmatrix}
     \int_{\xi=x_0}^{1} l_i(\xi) \, \diff \xi \\
     \int_{\xi=x_1 }^{1} l_i(\xi) \, \diff \xi \\
     \int_{\xi=x_2 }^{1} l_i(\xi) \, \diff \xi \\
     \vdots \\
     \int_{\xi=x_{N-1}}^{1} l_i(\xi) \, \diff \xi \\
     \int_{\xi=x_{N}}^{1} l_i(\xi) \, \diff \xi
   \end{bmatrix}.
 \end{equation}
The first element $\int_{\xi=x_0 }^{1} l_i(\xi) \, \diff \xi$ must be equal to zero, which is granted by the proposed construction/choice of the integration constant.
 
To conclude, we see that for given $i$ the computation of integrals $\int_{\xi=-1 }^{x_k} l_i(\xi) \, \diff \xi$ and $\int_{\xi=x_k }^{1} l_i(\xi) \, \diff \xi$ at \emph{all} collocation points $k=0, \dots, N$ requires us to compute two discrete cosine transforms, $\DCTI_{\text{Wolfram}}/\DCTI_{\text{Matlab}}$, and perform simple manipulations with vector elements namely $\INTEGRATE$, $\EXTEND$ and $\REDUCE$, see Algorithm~\ref{alg:discrete-dirac-integrals}. These integrals are the key component in the construction of the \emph{complete} $i$-th column of the Green matrix, see~\eqref{eq:93}. Furthermore, the computations for different values of $i$ are independent, and they can be run in parallel. Consequently, we see that the computation of Green matrix elements is numerically feasible even for large~$N$, that is for large number of collocation points/high polynomial degree. Clearly, the action of the discrete solution operator/Green matrix to the discrete right-hand side vector can be also implemented in a matrix free fashion.
 
We note that same integration method can not be used for the integrals  $\int_{\xi=-1 }^{x_k} \xi l_i(\xi) \, \diff \xi$. The integrand is in this case a polynomial of degree $N+1$, hence its exact representation in the Chebyshev--Fourier space would require us to know its values at $N+2$ collocation points/nodes. But since we have at our disposal only $N+1$ collocation points, and we can not easily expand the number of collocation points---for the different set of collocation points we would lose the property $l_i(x_j) = \kdelta{_{ij}}$. Thus the integrals $\int_{\xi=-1 }^{x_k} \xi l_i(\xi) \, \diff \xi$ must be computed by the formula \eqref{eq:34}.

\begin{algorithm}[!t]
  \caption{Compute $\int_{\xi=-1 }^{x_k} l_i(\xi) \, \diff \xi$ and $\int_{\xi=x_k }^{1} l_i(\xi) \, \diff \xi$, polynomial degree $N$, collocation points $\left\{x_k\right\}_{k=0}^N$}
  \label{alg:discrete-dirac-integrals}
\begin{algorithmic}[1]
  \Require $0 \leq i \leq N$ \Comment{We have $N+1$ Lagrange polynomials, numbering of $l_i$ starts from zero, \eqref{eq:20}.}
  \State
  $
  \vec{l}_i \gets
  \begin{bmatrix}
    \kdelta{_{(i+1)j}}
  \end{bmatrix}_{j=1}^{N+1}
  $
  \Comment{Node value representation of $l_i$, \eqref{eq:69}.}
  \State
  $\widehat{\vec{l}}_i
  \gets 
  \tensorschur{
    \left(\sqrt{\vec{c}}\right)^{-1}
  }
  {
    \DCTI_{\text{Matlab}}
    \left(
      \sqrt{\frac{2}{N}}
      \left(\sqrt{\vec{c}}\right)^{-1}
      \vec{l}_i
    \right)
  }
  $
  \Comment{Chebyshev--Fourier representation of $l_i$, \eqref{eq:71}.}
  \State $\widehat{\vec{l}}_{i, 2N} \gets \EXTEND \left( \widehat{\vec{l}}_i \right)$ \Comment{Extension of $\widehat{\vec{l}}_i$ by $N$ trailing zeros, \eqref{eq:74}.}
  \State
  $
  \widehat{\vec{L}}_{i, 2N}
  \gets
  \INTEGRATE
  \left(
    \widehat{\vec{l}}_{i, 2N}
  \right)
  $
  \Comment{Integral $\int l_i(\xi) \, \diff \xi$, Chebyshev--Fourier representation, \eqref{eq:51}.}
  \State
  $
   \vec{L}_{i, 2N}
  \gets
    \tensorschur{
     \sqrt{
         \vec{c}
       }
     }
     {
   \DCTI_{\text{Matlab}}
   \left(
     \sqrt{N}
     \tensorschur{
       \sqrt{
         \vec{c}
       }
     }
     {
       \widehat{\vec{L}}_{i, 2N}
     }
   \right)
 }
 $
 \Comment{Node value representation of $\int l_i(\xi) \, \diff \xi$, fine grid, \eqref{eq:77}.}
 \State
 $
 \vec{L}_{i, N}
 \gets
 \REDUCE
 \left(
   \vec{L}_{i, 2N}
 \right)
 $
 \Comment{Node values $\int l_i(\xi) \, \diff \xi$, original grid, \eqref{eq:84}.}
 
 \Comment{Fix the integration constants.}


 \State 
 $
 \vec{L}_{i, N}^{\uparrow}
 \gets
 \vec{L}_{i, N}
 -
 L_{i, N}(x_N)
 \vec{1}_{N+1}
 $
 \Comment{Node values, $\int_{\xi=x_N }^{x_k} l_i(\xi) \, \diff \xi = \int_{\xi=-1}^{x_k} l_i(\xi) \, \diff \xi$, \eqref{eq:86}.}
 \State
 $
   \vec{L}_{i, N}^{\downarrow}
   \gets 
   L_{i, N}(x_0)
   \vec{1}_{N+1}
   -
   \vec{L}_{i, N}
 $
 \Comment{Node values, $\int_{\xi=x_k }^{x_0} l_i(\xi) \, \diff \xi = \int_{\xi=x_k }^{1} l_i(\xi) \, \diff \xi$, \eqref{eq:88}.}

\end{algorithmic}
\end{algorithm}



\subsection{Evaluating Chebyshev polynomials at collocation points}
\label{sec:eval-chebysh-polyn}
Finally, we note that the values of Chebyshev polynomials at the Chebyshev--Gauss--Lobatto collocation points $\left\{ z_i \right\}_{i=0}^M$, which is a part of the computation~\eqref{eq:34}, can be found using the discrete cosine transform as well. For example, if we want to evaluate the Chebyshev polynomial $\ChebyshevT{k}$ at all points $\left\{ z_i \right\}_{i=0}^M$, then we first represent the polynomial in the Chebyshev--Fourier space,
\begin{equation}
  \label{eq:65}
  \widehat{\vec{T}}_k
  =
  \transpose{
  \begin{bmatrix}
    0 &
    \cdots &
    0 &
    1 &
    0 &
    \cdots &
    0
  \end{bmatrix}
},
\end{equation}
where the only nonzero element is at the $(k+1)$-th position. (Provided that we number the vector element from one. The dimension of the vector is $M+1$.) The discrete cosine transform of this vector
\begin{subequations}
  \label{eq:66}
  \begin{align}
    \label{eq:67}
    \vec{T}_k
   &=
   \DCTI_{\text{Wolfram}}
   \left(
     \sqrt{\frac{M}{2}}
     \tensorschur{
       \vec{c}
     }
     {
    \widehat{\vec{T}}_k
     }
   \right)
     ,
    \\
    \label{eq:90}
    \vec{T}_k
    &=
   \tensorschur{
     \sqrt{
         \vec{c}
       }
     }
     {
   \DCTI_{\text{Matlab}}
   \left(
     \sqrt{\frac{M}{2}}
     \tensorschur{
       \sqrt{
         \vec{c}
       }
     }
     {
    \widehat{\vec{T}}_k
     }
   \right)
   },
  \end{align}
\end{subequations}
gives us the vector with the sought node values,
\begin{equation}
  \label{eq:68}
  \vec{T}_k
  =
  \transpose{
  \begin{bmatrix}
    \ChebyshevT{k}(z_0) &
    \ChebyshevT{k}(z_1) &
    \ChebyshevT{k}(z_2) &
    \cdots &
    \ChebyshevT{k}(z_{M-1}) &
    \ChebyshevT{k}(z_{M})
  \end{bmatrix}
}.
\end{equation}
(We again need to choose the appropriate version of the discrete cosine transform implemented in the given software package.)

In particular, if we want to find the node values of the function
\begin{equation}
  \label{eq:91}
P_{i, N}
=_{\bydefinition}
\frac{\bweight_i}{2^{N+1}}
\left(
  \frac{\ChebyshevT{N+2}(x)}{N+2}
  -
  \frac{2 \ChebyshevT{N}(x)}{N}
  +
  \frac{\ChebyshevT{N-2}(x)}{N-2}
\right),  
\end{equation}
at all collocation points $\left\{ x_k \right\}_{k=0}^N$, see~\eqref{eq:34}, we can do it using Algorithm~\ref{alg:node-polynomial-int}. The algorithm is based on the observation that the Chebyshev--Fourier representation of $P_{i, N}$ can be written as a vector of dimension $2N+1$ with zero element everywhere except at positions $N-1$, $N+1$ and $N+3$. (Provided that the elements are numbered from one.) The node value representation at the fine set of collocation points~$\left\{ y_j \right\}_{j=0}^{2N}$ is then obtained by the discrete cosine transform~\eqref{eq:66} of this vector. Subsequently, one exploits the interlacing property between the original collocation points~$\left\{ x_k \right\}_{k=0}^N$ and the fine collocation points~$\left\{ y_j \right\}_{j=0}^{2N}$, see~\eqref{eq:81} and the detailed discussion therein, and one obtains the sought values of $P_{i, N}$ at original set of collocation points. We note that in implementation of the algorithm one can manually cancel the factor $2^N$ in the term~$\frac{\bweight_i}{2^{N+1}}$, see formulae~\eqref{eq:31} for the barycentric weights $\lambda_i$.

\begin{algorithm}[!t]
  \caption{
    Compute $P_{i, N}$, polynomial degree $N$, collocation points $\left\{x_k\right\}_{k=0}^N$
  }
  \label{alg:node-polynomial-int}
  \begin{algorithmic}[1]
    \Require $0 \leq i \leq N$ \Comment{We have $N+1$ Lagrange polynomials, numbering of $l_i$ starts from zero, \eqref{eq:20}.}
    \State
    $
    \widehat{\vec{P}}_{i,  2N} \gets
    \frac{\bweight_i}{2^{N+1}}
    \transpose{
      \begin{bmatrix}
        0 &
        \cdots &
        0 &
        \frac{1}{N-2} &
        0 &
        -\frac{2}{N} &
        0 &
        \frac{1}{N+2} &
        0 &
        \cdots \hspace{2\arraycolsep} 0 
      \end{bmatrix}
    }$
  \Comment{Chebyshev--Fourier representation or $P_{i, N}$, \eqref{eq:91}.}
  \State
  $
  \vec{P}_{i, 2N}
  \gets
  \tensorschur{
    \sqrt{
      \vec{c}
    }
  }
  {
    \DCTI_{\text{Matlab}}
    \left(
      \sqrt{N}
      \tensorschur{
        \sqrt{
          \vec{c}
        }
      }
      {
        \widehat{\vec{P}}_{i, 2N}
      }
    \right)
  }
  $
  \Comment{Node value representation, fine grid, \eqref{eq:90}.}
  \State
  $
  \vec{P}_{i, N}
  \gets
  \REDUCE \left( \vec{P}_{i, 2N} \right)
  $
  \Comment{Node values, original grid, \eqref{eq:85}.}

  \Comment{Fix the integration constants.}
  \State 
  $
  \vec{P}_{i, N}^{\uparrow}
  \gets
  \vec{P}_{i, N}
  -
  P_{i, N}(x_N)
  \vec{1}_{N+1}
  $
  \Comment{Node values, \eqref{eq:95}.}

  \State
  $
   \vec{P}_{i, N}^{\downarrow}
   \gets 
   P_{i, N}(x_0)
   \vec{1}_{N+1}
   -
   \vec{P}_{i, N}
   $
   \Comment{Node values, \eqref{eq:94}.}
 \end{algorithmic}
\end{algorithm}

\subsection{Algorithm}
\label{sec:algorithm}

The algorithm for the computation of Green matrix elements is summarised in Algorithm~\ref{alg:green-matrix-elements}. The implementation in Wolfram Mathematica is almost verbatim transcription of the pseudocode in Algorithm~\ref{alg:green-matrix-elements}, the only substantial difference is that we use vectorisation capabilities of the programming language. On the other hand, the implementation in Matlab is aimed at speed, and in addition to vectorisation capabilities of the programming language it also exploits the centrosymmetry of the Green matrix and it uses some other optimisations, see the comments in the source code.

Finally, we note that that application of Green matrix to a given vector can be coded in a \emph{matrix-free form}, and we implement this operation in our proof-of-concept Matlab code as well. In this regard the situation is the same as for the spectral differentiation. The application of the second order spectral differentiation matrix to a vector also has a matrix-free variant, see~\cite[Section 3.4.3]{shen.j.tang.t.ea:spectral}. Moreover, the matrix-free differentiation is in this case implemented also via the backward/forward discrete cosine transformation and the differentiation itself takes place in the Chebyshev--Fourier space by simple manipulations with Chebyshev--Fourier coefficients. The algorithmic complexity of the matrix-free spectral differentiation is known to be $\bigo{N \log_2 N}$, see for example~\cite[Section 3.4.3]{shen.j.tang.t.ea:spectral}, and the matrix-free application of discrete Green operator is of the same complexity. Indeed, the matrix-free version of Algorithm~\ref{alg:green-matrix-elements} also utilises---as the spectral matrix-free differentiation algorithm---forward/backward cosine transformation and a simple manipulation in Chebyshev--Fourier space.

\begin{algorithm}[!t]
  \caption{
    Compute Green matrix $\GMatrix{N}$,  polynomial degree $N$, collocation points $\left\{x_k\right\}_{k=0}^N$
  }
  \label{alg:green-matrix-elements}
  \begin{algorithmic}[1]
    \State Get barycentric weights $\vec{\lambda}$,  $\left\{ \lambda_k \right\}_{k=0}^N$ \Comment{Barycentric weights, \eqref{eq:15}}
    \State Get Chebyshev--Gauss--Lobatto collocation points $\vec{x}$, $\left\{ x_k \right\}_{k=0}^N$ \Comment{Collocation points, \eqref{eq:12}.}

    \For{$0 \leq i \leq N$}     \Comment{Get $i$-th column of Green matrix, columns are numbered from zero.}
    \State Get $\vec{L}_{i, N}^{\downarrow}$, $\vec{L}_{i, N}^{\uparrow}$ \Comment{Use Algorithm~\ref{alg:discrete-dirac-integrals}.}
    \State Get $\vec{P}_{i, N}^{\downarrow}$, $\vec{P}_{i, N}^{\uparrow}$ \Comment{Use Algorithm~\ref{alg:node-polynomial-int}.}
    \State
    $
    \vec{G}_{i}\gets 
    \frac{1}{2}
    \tensorschur{
    \left(
      \vec{x} + \vec{1}_{N+1}
    \right)
  }
  {
    \left[
      \vec{P}_{i, N}^{\downarrow}
      +
      \left(x_i - 1\right)
      \vec{L}_{i, N}^{\downarrow}
    \right]
  }
    +
    \frac{1}{2}
    \tensorschur{
    \left(
      \vec{x} - \vec{1}_{N+1}
    \right)
  }
  {
    \left[
      \vec{P}_{i, N}^{\uparrow}
      +
      \left(x_i + 1\right)
      \vec{L}_{i, N}^{\uparrow}
    \right]
  }
    $
    \Comment{Evaluate~\eqref{eq:93}.}
    \EndFor
    \State
    $
    \GMatrix{N}
    \gets
    \begin{bmatrix}
      \vec{G}_{0} &
      \vec{G}_{1} &
      \cdots &
      \vec{G}_{N}
    \end{bmatrix}
    $
    \Comment{Collect all columns.}
  \end{algorithmic}
\end{algorithm}

\subsection{Green matrix as an inverse matrix to the second order spectral differentiation matrix}
\label{sec:green-matrix-as}
Clearly, the second order spectral differentiation matrix $\DMatrix{N}^2$ as such is singular since it has a non-trivial kernel corresponding to all affine functions. Thus the invertibility of $\DMatrix{N}^2$ must be understood in the context of solution of the discretised problem, that is including the boundary conditions. Below we introduce three approaches to the concept of ``inverse'' to the second order spectral differentiation matrix. 

\subsubsection{Left inverse}
\label{sec:left-inverse}
If we decide to implement the boundary conditions by stripping the first and the last row and the first and the last column from the spectral differentiation matrix, see, for example, \cite[Chapter 7]{trefethen.ln:spectral} for details, then the matrix $\DMatrix{N}^2$ is reduced to matrix $\DDMatrixTilde{N}$, and the corresponding algebraic problem reads
\begin{equation}
  \label{eq:106}
  \DDMatrixTilde{N} \widetilde{\vec{y}}= \widetilde{\vec{f}},
\end{equation}
where the vectors $\widetilde{\vec{y}}$ and $\widetilde{\vec{f}}$ are vectors of node values of functions $y$ and $f$ at the \emph{inner} collocation points $\left\{x_k \right\}_{k=1}^{N-1}$, see Section~\ref{sec:standard-approach}. For polynomial order $N$ the problem~\eqref{eq:106} is an algebraic problem with a matrix of dimension $(N-1) \times (N-1)$.

On the other hand the Green matrix $\GMatrix{N}$ has the boundary conditions already built in, and for polynomial order $N$ it is a matrix of dimension $(N+1) \times (N+1)$ with the zero first and last row. The node values of solution $y$ at \emph{all} collocation points $\left\{x_k \right\}_{k=0}^{N}$ are given by the application of $\GMatrix{N}$ to the node values of $f$ at \emph{all} collocation points as well,
\begin{equation}
  \label{eq:107}
  \vec{y} = \GMatrix{N} \vec{f}.
\end{equation}

Clearly, there is an incompatibility in between the dimensions of $\DDMatrixTilde{N}$ and $\GMatrix{N}$, hence one matrix can hardly be the inverse of the other. The invertibility property between the matrices however holds if we restrict ourselves to the inner collocations points. At inner collocation points we have
\begin{equation}
  \label{eq:109}
  \widetilde{\GMatrix{N} \DMatrix{N}^2} = \identity_{(N-1) \times (N-1)},
\end{equation}
where $\identity_{(N-1) \times (N-1)}$ denotes the identity matrix of dimension $(N-1) \times (N-1)$. This means that if we first take the product of the Green matrix $\GMatrix{N}$ and the \emph{full} second order spectral differentiation matrix~$\DMatrix{N}^2$, and then we strip the so-obtained matrix  of its first and last row and its first and last column, then we get the identity matrix of dimension $(N-1) \times (N-1)$. In this sense the Green matrix $\GMatrix{N}$ the left inverse of the second order spectral differentiation matrix~$\DMatrix{N}^2$.

\subsubsection{Right inverse}
\label{sec:right-inverse}
The right inverse should be the counterpart to the continuous level calculation~\eqref{eq:17}. However, we again can not investigate directly the product $\DMatrix{N}^2\GMatrix{N}$. If we apply the Green operator~\eqref{eq:111} to a polynomial of degree at most $N$, or in other words if we solve~\eqref{eq:1} for the right-hand side being a polynomial of degree at most $N$, then we get the solution that is of polynomial degree at most $N+2$. The application of spectral second order differentiation matrix~$\DMatrix{N}^2$ is however exact only for polynomials of order $N$, hence we can not expect $\DMatrix{N}^2\GMatrix{N} = \identity_{(N+1) \times (N+1)}$. We must adjust the polynomial degree accordingly.

We start with a polynomial $p_{N-2}$ of degree at most $N-2$ interpolated at Chebyshev--Gauss--Lobatto points $\left\{z_j\right\}_{j=0}^{N-2}$. At this set of collocation points the polynomial is \emph{exactly} represented by the vector $\vec{p}_{N-2}$. (Meaning that there is one-to-one correspondence between polynomials of degree at most $N-2$ and their node values at $\left\{z_j\right\}_{j=0}^{N-2}$.) If we reinterpolate the polynomial to a finer set of Chebyshev--Gauss--Lobatto points $\left\{x_k \right\}_{k=0}^{N}$, then we still have an \emph{exact} representation of this polynomial $p_{N-2}$, but with a longer vector $\vec{p}_{N}$. We denote this reinterpolation operation as
\begin{equation}
  \label{eq:110}
  \vec{p}_N = \tensorq{R}_{N-2 \to N} \vec{p}_{N-2},
\end{equation}
and we recall that software tools such as \texttt{chebfun} package have an infrastructure for this operation, see \texttt{barymat} function in~\texttt{chebfun}.

The reinterpolated polynomial can be acted upon by the Green matrix~$\GMatrix{N}$. The action of Green operator~$\GOperator$ to a polynomial~$p_{N-2}$ of degree~$N-2$ produces a polynomial of degree $N$, thus on the discrete level we still have \emph{exact} representation of $\GMatrix{N}\tensorq{R}_{N-2 \to N} \vec{p}_{N-2}$ on the finer set of collocation points. (The use of finer set of collocation points than is necessary to represent $\vec{p}_{N-2}$ compensates the fact that the application of Green operator increases the polynomial order by two.) The action of the second order derivative operator, that is the application of the second order spectral differentiation matrix $\DMatrix{N}^2$, is exact for degree $N$ polynomials. Consequently, the discrete level operation
\begin{equation}
  \label{eq:112}
  \DMatrix{N}^2 \GMatrix{N}\tensorq{R}_{N-2 \to N} \vec{p}_{N-2}
\end{equation}
is the exact representation of $\ddd{}{x} \GOperator(p_{N-2})$, and the result of this operation is a polynomial of degree $N-2$. This polynomial can be reinterpolated to the coarse set of collocation points by the reinterpolation operator $\tensorq{R}_{N \to N-2}$. The sequence of operations
\begin{equation}
  \label{eq:123}
  \tensorq{R}_{N \to N-2} \DMatrix{N}^2 \GMatrix{N} \tensorq{R}_{N-2 \to N}
\end{equation}
defines a \emph{square} matrix of dimension $(N-2)\times(N-2)$, and this matrix is an exact representation of the operator $\ddd{}{x} \GOperator$ at the set of polynomials of degree at most $N-2$. Thus we have
\begin{equation}
  \label{eq:124}
  \tensorq{R}_{N \to N-2} \DMatrix{N}^2 \GMatrix{N} \tensorq{R}_{N-2 \to N} = \identity_{(N-2)\times(N-2)}.
\end{equation}
In this sense the Green matrix $\GMatrix{N}$ is the right inverse of the second order spectral differentiation matrix~$\DMatrix{N}^2$.

\subsubsection{Left and right inverse in terms of reduction/extension to/from inner collocation points}
\label{sec:left-right-inverse}

Yet another option how to interpret $\GMatrix{N}$ as the inverse to $\DMatrix{N}^2$ and vice versa is to \emph{embed the boundary conditions into the differentiation matrix~$\DMatrix{N}^2$}. See~\cite{aurentz.jl.trefethen.ln:block} for a similar construction.

If we apply the second derivative matrix $\DMatrix{N}^2$ to a polynomial  $p_{N}$ of degree $N$ represented by its values $\vec{p}_{N}$ at all $N+1$ Chebyshev collocation points, we get the polynomial
\begin{equation}
  \label{eq:170}
  p_{N-2} =_{\bydefinition} \ddd{}{x} p_{N}
\end{equation}
of degree $N-2$ represented by its node values $\vec{p}_{N-2}$,
\begin{equation}
  \label{eq:171}
  \vec{p}_{N-2} = \DMatrix{N}^2 \vec{p}_{N},
\end{equation}
at all $N+1$ Chebyshev collocation points. However, for the exact characterisation of this $N-2$ degree polynomial, we need to know only $N-1$ values of this polynomial at $N-1$ collocation points. In this sense the $N-2$ degree polynomial $p_{N-2}$ is overrepresented provided that the polynomial interpolation uses node values at \emph{complete collection} of $N+1$ Chebyshev collocation points. Thus we may decide to represent the polynomial $p_{N-2}$ only by its values at the $N-1$ \emph{inner} Chebyshev collocations points.

The reduction operation is represented by the matrix $\tensorq{P}_{N \to N-2}$ that just crops the vector $\vec{p}_{N-2}$ out of its first and last element, that is
\begin{equation}
  \label{eq:163}
  \tensorq{P}_{N \to N-2}
  =
  \begin{bmatrix}
    \vec{0}_{N-1} & \identity_{(N-2)\times(N-2)} & \vec{0}_{N-1}
  \end{bmatrix}
  =
  \begin{bmatrix}
    \vspace{-0.3em}
    \includegraphics[scale=0.5]{figures-matrix-structure/p-matrix-structure}                    
  \end{bmatrix}_{\left(N-1\right) \times \left(N+1\right)}
  ,
\end{equation}
where $\identity_{(N-2)\times(N-2)}$ is the square identity matrix of dimension $(N-2)\times(N-2)$ and~$\vec{0}_{N-1}$ denotes the zero column vector of dimension $N-1$. The simple projection matrix $\tensorq{P}_{N \to N-2}$ is a rectangular matrix of dimension $\left(N-1\right) \times \left(N+1\right)$, and it reduces vector $\vec{p}_{N-2}$ of dimension $N+1$ to the vector of dimension $N-1$. This vector still represents the same polynomial and the representation is exact, nothing has been lost. Note that unlike in Section~\eqref{sec:right-inverse} \emph{we are not reinterpolating to a new set of collocation points}, we just use only \emph{inner collocation points} out of the original ones. (The use of inner collocation points also makes the present approach different from~\cite{aurentz.jl.trefethen.ln:block}. \cite{aurentz.jl.trefethen.ln:block} who also use reduced vectors to represent the same polynomial, but their approach requires \emph{reinterpolation} to the coarser set of \emph{new} collocation points. The reinterpolation is avoided in the present approach.)  The $\left(N-1\right) \times \left(N+1\right)$ rectangular matrix obtained by the multiplication
\begin{equation}
  \label{eq:164}
  \tensorq{P}_{N \to N-2} \DMatrix{N}^2
\end{equation}
thus evaluates the second derivative at the $N-1$ inner collocation points. This operation is guaranteed to be an exact representation of the second derivative for polynomials of order at most $N$. Since $\tensorq{P}_{N \to N-2} \DMatrix{N}^2$ is rectangular matrix, we can add two rows to make it a square matrix.

These two additional rows are utilised for the boundary conditions. The square matrix $\DDMatrixBC{N}$ defined as 
\begin{equation}
  \label{eq:165}
  \DDMatrixBC{N}
  =_{\bydefinition}
  \begin{bmatrix}
    \begin{matrix} 1 & \transpose{\vec{0}_{N}} \end{matrix} \\
    \tensorq{P}_{N \to N-2} \DMatrix{N}^2 \\
    \begin{matrix} \transpose{\vec{0}_{N}} & 1 \end{matrix}
  \end{bmatrix}
  =
  \begin{bmatrix}
    \vspace{-0.3em}
    \includegraphics[scale=0.5]{figures-matrix-structure/ddmatrixbc-matrix-structure}                    
  \end{bmatrix}_{(N+1) \times (N+1)}
  ,
\end{equation}
where $\transpose{\vec{0}_{N}}$ denotes the zero row vector of dimension $N$, represents both the second derivative operator $\ddd{}{x}$ and evaluation operator at the end points of $[-1, 1]$. In other words, $\DDMatrixBC{N} \vec{p}_{N}$ returns a column vector where the first element is just a copy of the first element of the input vector $\vec{p}_{N}$ (the value of $p_N$ at collocation point $x_0$), while the last element is just a copy of the first element of the input vector $\vec{p}_{N}$ (the value of $p_N$ at collocation point $x_N$). The remaining elements represent the values of the second derivative at the inner collocation points.

An accompanying operator to the \emph{reduction} operator $\tensorq{P}_{N \to N-2}$ is the \emph{extension} operator~$\tensorq{E}_{N-2 \to N}$. The extension operator takes~$N-1$ values of polynomial $p_{N-2}$ at the inner collocation points and maps them to the values of (the same polynomial)~$p_{N-2}$ at the complete collection of~$N+1$ collocation points. Having the \emph{values} of polynomial $p_{N-2}$ of degree $N-2$ at the $N-1$ inner collocation points, we can uniquely construct the polynomial $p_{N-2}$ using the Lagrange interpolant, and evaluate the interpolant at the endpoints $x_0$ and $x_{N}$ of the complete collection of collocation points. (In evaluating the interpolation polynomial we naturally use the barycentric interpolation formula, in particular in \texttt{chebfun} we use \texttt{bary} function.) We denote the matrix representation of this extension operator as~$\tensorq{E}_{N-2 \to N}$,
\begin{equation}
  \label{eq:166}
  \tensorq{E}_{N-2 \to N}
  =_{\bydefinition}
  \begin{bmatrix}
    \vec{I}_{N-1, 0} \\
    \identity_{(N-2)\times(N-2)} \\
    \vec{I}_{N-1, N}
  \end{bmatrix}
  =
  \begin{bmatrix}
    \vspace{-0.3em}
    \includegraphics[scale=0.5]{figures-matrix-structure/e-matrix-structure}                    
  \end{bmatrix}_{\left(N+1\right) \times \left(N-1\right)}
  ,
\end{equation}
where $\vec{I}_{N-1, 0}$ and $\vec{I}_{N-1, N}$ represent row vectors that correspond to the construction of Lagrange interpolant and its evaluation at $x_0$ and $x_N$ respectively.

Having defined the extension operator/matrix, we can introduce the modified Green matrix as
\begin{equation}
  \label{eq:167}
  \GMatrixBC{N} =_{\bydefinition}
  \begin{bmatrix}
    \frac{1}{2} \left( x_0 + \vec{x}  \right)
    &
    \GMatrix{N}\tensorq{E}_{N-2 \to N}
    &
    -\frac{1}{2} \left( x_N + \vec{x}\right)
  \end{bmatrix}
  =
  \begin{bmatrix}
    \vspace{-0.3em}
    \includegraphics[scale=0.5]{figures-matrix-structure/gmatrixbc-matrix-structure}                    
  \end{bmatrix}_{(N+1) \times (N+1)}
  ,
\end{equation}
where the first and the last column represent \emph{linear} functions that have two properties. First, they vanish at $x_0$ and $x_N$ respectively---the last/first elements in the corresponding vector representation are equal to zero. Second, they take value one at $x_N$ and $x_0$ respectively---the first/last elements in the corresponding vector representation are equal to one. Note that since these vectors represent linear functions, they vanish upon the action of the second derivative operator/matrix.

The block structure for the product of the matrices $\DDMatrixBC{N}$ and $\GMatrixBC{N}$ is
\begin{equation}
  \label{eq:172}
  \DDMatrixBC{N} \GMatrixBC{N}
  =
  \begin{bmatrix}
    \vspace{-0.3em}
    \includegraphics[scale=0.5]{figures-matrix-structure/ddmatrixbc-matrix-structure}                    
  \end{bmatrix}_{(N+1) \times (N+1)}
  \begin{bmatrix}
    \vspace{-0.3em}
    \includegraphics[scale=0.5]{figures-matrix-structure/gmatrixbc-matrix-structure}                    
  \end{bmatrix}_{(N+1) \times (N+1)}
  ,
\end{equation}
and for $\DDMatrixBC{N}$ and $\GMatrixBC{N}$ we have the desired property
\begin{equation}
  \label{eq:168}
  \DDMatrixBC{N} \GMatrixBC{N}
  =
  \identity_{\left(N+1\right) \times \left(N+1\right)}
\end{equation}
and simultaneously 
\begin{equation}
  \label{eq:169}
  \GMatrixBC{N} \DDMatrixBC{N}  = \identity_{\left(N+1\right) \times \left(N+1\right)}.
\end{equation}

\section{Conclusion}
\label{sec:conclusion}

The Green matrix for the problem~\eqref{eq:1} is a discrete counterpart of the Green operator/solution operator for the problem~\eqref{eq:1}. We have proposed a fast algorithm for the construction of the Green matrix, and we have implemented the algorithm in Matlab and Wolfram Mathematica. The Matlab implementation provides interoperability with the \texttt{chebfun} package, which is a standard Matlab package used in, among others, implementation of spectral collocation method. The Wolfram Mathematica implementation provides interoperability with symbolic analytical tools available in this software. The \emph{proof-of-concept} implementation in Matlab shows that the algorithm is indeed applicable for a large number of collocation points/high polynomial degree, and that the action of the discrete Green operator can be implemented in a matrix-free fashion. Furthermore, the computational time needed to construct the Green matrix is comparable to the time needed to construct the second order spectral differentiation matrix via \texttt{diffmat} function in \texttt{chebfun} package.

\appendix

\section{Consistent scalar product and centrosymmetry of second order spectral differentiation matrix}
\label{sec:cons-scal-prod}

\subsection{Consistent scalar product}
\label{sec:cons-scal-prod-1}
On the continuous level we equip the appropriate function space with the standard scalar product
\begin{equation}
  \label{eq:142}
  \fascalalt{u}{v} =_{\bydefinition} \int_{x=-1}^1 u(x) v(x) \, \diff x, 
\end{equation}
and we want to introduce a discrete level scalar product that inherits the properties of~\eqref{eq:142}. In particular, given two polynomials $p$, $q$ of degree at most $N$, we want to calculate the exact value of $\fascalalt{p}{q}$ using the values of $p$ and $q$ at the Chebyshev--Gauss--Lobatto points collocation points $\left\{ x_j \right\}_{j=0}^N$.

The standard \emph{Clenshaw--Curtis quadrature rule} for a generic function $f$ at this set of collocation points reads
\begin{equation}
  \label{eq:152}
  \int_{x=-1}^1 f(x) \, \diff x
  \approx
  \sum_{i=0}^N f(x_i) \qweight_i,
\end{equation}
where~$\left\{ w_{N, i} \right\}_{i=0}^N$ are the quadrature weights, see~\cite{trefethen.ln:approximation} or~\cite{mason.jc.handscomb.dc:chebyshev}, for the collocation points~$\left\{ x_j \right\}_{j=0}^N$. (Note that we indeed use Clenshaw--Curtis quadrature, not the Chebyshev--Gauss--Lobatto quadrature.) The quadrature~\eqref{eq:152} is exact provided that $f$ is a polynomial of degree at most $N$. In the scalar product~$\fascalalt{p}{q}$ we are however multiplying two polynomials of degree at most $N$, which can lead to a polynomial of degree at most $2N$, and we can not, in general, expect the exactness of the product in the sense
\begin{equation}
  \label{eq:21}
  \fascalalt{p}{q} = \int_{x=-1}^1 p(x) q(x) \, \diff x = \sum_{i=0}^N p(x_i) q(x_i) \qweight_{N, i}.
\end{equation}
The exactness can be reached provided that we reinterpolate, at the discrete level, to a finer set of collocation points. This can be done using the reinterpolation operator~$\tensorq{R}_{N \to 2N}$ introduced in Section~\ref{sec:green-matrix-as}. Namely, if we for the polynomial~$p$ use vector~$\vec{p}$ as its node value representation at collocation points $\left\{ x_j \right\}_{j=0}^N$, then its representative at the collocation points~$\left\{ y_k \right\}_{k=0}^{2N}$ is
\begin{equation}
  \label{eq:153}
  \tensorq{R}_{N \to 2N}  \vec{p} 
\end{equation}
or in components $\tensor{\left(\tensorq{R}_{N \to 2N}  \vec{p} \right)}{_{i}} = \sum_{j=0}^N\tensor{\left(\tensorc{R_{N \to 2N}}\right)}{_{ij}}p_j$. (This operation converts the vector of dimension $N+1$ to the vector of dimension $2N+1$, but this vector is still representing the same function, the polynomial $p$.) We do the same for the polynomial $q$, and in the integral~\eqref{eq:142} we use the Clenshaw--Curtis quadrature at the finer set of collocation points which leads to
\begin{multline}
  \label{eq:154}
  \fascalalt{p}{q} = \int_{x=-1}^1 p(x) q(x) \, \diff x
  =
  \sum_{i,m=0}^{2N}
  \kdelta{_{im}}
  \left(
    \sum_{j=0}^N\tensor{\left(\tensorc{R_{N \to 2N}}\right)}{_{ij}}p_j
  \right)
  \left(
    \sum_{k=0}^N\tensor{\left(\tensorc{R_{N \to 2N}}\right)}{_{mk}}q_k
  \right)
  \qweight_{2N, i}
  =
  \vectordot{
    \left(
    \tensorq{W}_{2N}
    \tensorq{R}_{N \to 2N}
    \vec{p}
    \right)
  }
  {
    \left(
      \tensorq{R}_{N \to 2N}
      \vec{q}
    \right)
  }
  \\
  =
  \vectordot{
    \left(
      \transpose{\tensorq{R}_{N \to 2N}}
      \tensorq{W}_{2N}
      \tensorq{R}_{N \to 2N}
      \vec{p}
    \right)
  }
  {
    \vec{q}
  }
  =
  \vectordot{\tensorq{S}_{2N}\vec{p}}{\vec{q}}
  .
\end{multline}
(The Clenshaw--Curtis quadrature now deals with the product $pq$ which is still a polynomial of degree at most $2N$, but represented at the finer set of collocation points wherein the quadrature is exact.) In~\eqref{eq:154} we use the symbol
\begin{equation}
  \label{eq:161}
  \vectordot{\vec{a}}{\vec{b}} = \sum_{i,j=0}^{N} \kdelta{_{ij}} \tensor{\tensorc{a}}{_i} \tensor{\tensorc{b}}{_j}
\end{equation}
for the standard scalar product on $\R^{N+1}$, and we introduce matrices
\begin{subequations}
  \label{eq:155}
  \begin{align}
    \label{eq:156}
    \tensorq{W}_{2N} &=_{\bydefinition}
                       \begin{bmatrix}
                         \qweight_{2N, 0} & 0 & 0 & \cdots & 0 \\
                         0  & \qweight_{2N, 1} & 0 & \cdots & 0 \\
                         0 & 0 & \qweight_{2N, 2} & \cdots & 0 \\
                         \vdots & \vdots & \vdots & \ddots & \vdots \\
                         0 & 0 & 0 & \cdots & \qweight_{2N, 2N} \\
                       \end{bmatrix}
    ,
    \\
    \label{eq:157}
    \tensorq{S}_{2N} &=_{\bydefinition}
                  \transpose{\tensorq{R}_{N \to 2N}}
                  \tensorq{W}_{2N}
                  \tensorq{R}_{N \to 2N}.
  \end{align}
\end{subequations}
The symbol $\tensorq{W}_{2N}$ denotes a diagonal matrix with the diagonal occupied by the Clenshaw--Curtis quadrature weights $\left\{ \qweight_{2N, i} \right\}_{i=0}^{2N}$ at the \emph{finer} set of collocation points~$\left\{ y_k \right\}_{k=0}^{2N}$, and $\tensorq{S}_{2N}$ denotes the \emph{symmetric positive definite matrix} of dimension $(N+1) \times (N+1)$ that induces a new scalar product on the discrete level. (We recall that the Clenshaw--Curtis quadrature weights are known to be positivive, see~\cite{imhof.jp:on}.) Thus if we define the new scalar product
\begin{equation}
  \label{eq:158}
  \vectordotalt{\vec{p}}{\vec{q}}{\tensorq{S}_{2N}} =_{\bydefinition} \vectordot{\tensorq{S}_{2N}\vec{p}}{\vec{q}},
\end{equation}
then we see that the new scalar product has the desired property. If $p$ and $q$ are polynomials of order at most $N$ represented by their node values vectors $\vec{p}$ and $\vec{q}$ at the coarse set of collocation points~$\left\{ x_j \right\}_{j=0}^N$, then we have the equality
\begin{equation}
  \label{eq:159}
  \fascalalt{p}{q} = \int_{x=-1}^1 p(x) q(x) \, \diff x = \vectordotalt{\vec{p}}{\vec{q}}{\tensorq{S}_{2N}}.
\end{equation}
In this sense the new discrete scalar product~\eqref{eq:158} is the natural discretisation of the continuous level scalar product~\eqref{eq:142}. Note that we use reinterpolation to the grid wherein the operations of interest are exact for polynomial of order at most $2N$---this is the minimal requirement. But we could have used reinterpolation with any $K$ for which we have  $K \geq 2N$.

\subsection{Second order spectral differentiation matrix}
\label{sec:second-order-spectr}
The second order derivative operator $\ddd{}{x}$ is, on the continuous level, a Hermitian (symmetric) operator in the sense that
\begin{equation}
  \label{eq:160}
  \fascalalt{\ddd{f}{x}}{g}
  =
  \fascalalt{f}{\ddd{g}{x}}
  .
\end{equation}
(We work with zero Dirichlet boundary conditions.) Consequently, one can expect that the discrete counterpart of $\ddd{}{x}$ is a symmetric matrix. This is true for example for the discretisation based on finite differences. But the spectral differentiation matrix~$\DMatrix{N}^2$, see Section~\ref{sec:standard-approach} is \emph{not symmetric}. This lack of symmetry might be puzzling, but it is actually fine.

The symmetry in the standard Euclidean scalar product~\eqref{eq:161} is in fact of no interest since the standard Euclidean scalar product in $\R^{N+1}$ is not a proper discrete counterpart of the continuous level scalar product~\eqref{eq:142} for which we have~\eqref{eq:160}. But if we use the consistent scalar product~\eqref{eq:158}, we get the desired property
\begin{equation}
  \label{eq:162}
  \vectordotalt{\left( \DMatrix{N}^2 \vec{p} \right)}{\vec{q}}{\tensorq{S}_{2N}}
  =
  \vectordotalt{\vec{p}}{\left( \DMatrix{N}^2 \vec{q} \right)}{\tensorq{S}_{2N}}
  ,
\end{equation}
that holds for all $p$, $q$ polynomials of order at most $N$. Identity~\eqref{eq:162} is the sought counterpart of~\eqref{eq:160} on the discrete level. Thus we can conclude that the discrete second order spectral differentiation matrix \emph{is symmetric} in the right \emph{scalar product}.

The proof of~\eqref{eq:162} is straightforward. It suffices to realise that $\DMatrix{N}^2 \vec{p}$ and $\DMatrix{N}^2 \vec{q}$  are the exact node value representatives of $\ddd{p}{x}$ and $\ddd{q}{x}$ respectively, and that the second derivatives are polynomials of order at most $N-2$. Consequently, one can use~\eqref{eq:159} and the proof is complete.


\bibliographystyle{chicago}
\bibliography{vit-prusa}

\end{document}